\documentclass[10pt]{amsart}

\newcommand{\dash}{\nobreakdash-\hspace{0pt}}
\newcommand{\YD}[1]{_{#1}^{#1}\mathcal{YD}}
\newcommand{\VV}[1]{V_{#1}}
\newcommand{\Nic}{\mathcal{B}}
\newcommand{\nilcoxeter}{N_W}
\newcommand{\field}{\mathbb{C}}
\newcommand{\Braid}{\mathbb{B}}
\newcommand{\Symm}{\mathbb{S}}

\newcommand{\h}{\mathfrak{h}}
\newcommand{\rderiv}[1]{{\overleftarrow D}_{#1}}
\newcommand{\mult}{\cdot} 
\newcommand{\demazure}[1]{\partial_{#1}}
\newcommand{\rdemazure}[1]{{{\overleftarrow\partial\hspace{-0.25em}}_{#1}}}
\newcommand{\len}{\ell}

\newcommand{\Hom}{\mathrm{Hom}}
\newcommand{\qd}{{\mathit{quad}}}
\newcommand{\gam}[1]{\gamma_{#1}}

\newcommand{\supp}{\mathit{supp}\ }
\DeclareMathOperator{\id}{id}
\DeclareMathOperator{\End}{End}
\DeclareMathOperator{\im}{im}

\newenvironment{xproposition}{\protect\subsection{Proposition}\it}{}
\newenvironment{xlemma}{\protect\subsection{Lemma}\it}{}   
\newenvironment{xtheorem}{\protect\subsection{Theorem}\it}{}   
\newenvironment{xcriterion}{\protect\subsection{Criterion}\it}{}   
\newenvironment{xcorollary}{\protect\subsection{Corollary}\it}{}
\newenvironment{xconjecture}{\protect\subsection{Conjecture}\it}{}
\newenvironment{xremark}{\protect\subsection{Remark}}{}

\advance\oddsidemargin by-0.3in
\advance\evensidemargin by-0.3in
\advance\textwidth by 0.6in

\begin{document}
\author{Yuri Bazlov}
\address{School of Mathematical Sciences,
Queen Mary, University of London,
Mile End Road, London E1 4NS, UK}
\email{y.bazlov@qmul.ac.uk}
\begin{title}
{Nichols-Woronowicz algebra model for Schubert calculus on Coxeter groups}
\end{title}

\begin{abstract}
We realise the cohomology ring of a flag manifold, more generally 
the coinvariant algebra of an arbitrary finite Coxeter group 
$W$, as a commutative subalgebra of 
a certain Nichols-Woronowicz algebra in the Yetter-Drinfeld category
over $W$. This gives a braided Hopf algebra version of the
corresponding Schubert calculus.
The nilCoxeter algebra and its action on the coinvariant algebra by
divided difference operators are also realised 
in the Nichols-Woronowicz algebra. 
We discuss the relationship between
Fomin-Kirillov quadratic algebras, Kirillov-Maeno bracket algebras and
our construction.   
\end{abstract}
\subjclass[2000]{Primary 20G42, secondary 20F55, 14M15}
\maketitle
\setcounter{tocdepth}{1}
\tableofcontents
\pagestyle{myheadings}
\markboth{Y. Bazlov}{Nichols-Woronowicz algebra model for Schubert calculus}

\section{Introduction}

Few years ago, a new approach to the cohomology rings of the flag manifolds 
was developed by Fomin and Kirillov. In \cite{FK} they
introduced  a family of noncommutative graded algebras
$\mathcal{E}_n$, 
defined only by quadratic relations, 
such that the cohomology ring of the manifold 
$\mathit{Fl}_n$ 
of complete flags in $\field^n$ is realised 
as a graded commutative subalgebra of
$\mathcal{E}_n$. 
The Fomin\dash Kirillov algebras have
many other interesting properties, notably they are braided Hopf
algebras over the symmetric groups $\Symm_n$.  

The manifolds $\mathit{Fl}_n$ naturally correspond to the symmetric
groups $\Symm_n$ and have a version 
where $\Symm_n$ is replaced by any crystallographic Coxeter group $W$:
if $G$ is a semisimple Lie group whose Weyl group is $W$, and $B$
is the Borel subgroup of $G$, the homogeneous space $G/B$ is
the flag manifold of $G$. The cohomology ring of $G/B$ was shown by
Borel to be isomorphic to the coinvariant algebra $S_W$ of $W$. Hence
the cohomology of the flag manifold has a description purely in terms
of the invariant theory of 
$W$ --- we refer to this description as a Schubert calculus over $W$. 
If one looks only at the algebraic side of the picture,
$W$ does not need be crystallographic (see the book \cite{H} of
Hiller). Thus, any finite Coxeter group $W$    
admits a Schubert calculus, although not necessarily coming from the
cohomology of a geometric object.
 
It is then natural to try and extend the Fomin\dash Kirillov
construction to arbitrary Coxeter groups. 
Recently, Kirillov and Maeno 
suggested a generalisation of
$\mathcal{E}_n$, where the symmetric group $\Symm_n$ is replaced by a
finite Coxeter group $W$ with a set $S$ of Coxeter generators. 
The bracket algebras 
$\mathit{BE}(W,S)$, defined in \cite{KM}, are in general not
quadratic,
but the relations in $\mathit{BE}(W,S)$ are still given explicitly in terms of
the root system of $W$. 
However, the conjecture
that $\mathit{BE}(W,S)$ contains a copy of $S_W$ 
was verified in \cite{KM} only for classical crystallographic groups
$W$ and for $W$ of type $G_2$. 

In the present paper, we suggest a new and uniform construction for the
coinvariant algebra $S_W$ of an arbitrary Coxeter group $W$. 
We realise $S_W$ as a graded commutative subalgebra in a 
Nichols\dash Woronowicz  algebra $\Nic_W$, which itself is a braided
Hopf algebra over the group $W$ with a number of additional properties. 
Moreover, the so\dash called nilCoxeter algebra $\nilcoxeter$ 
and its important representation on $S_W$ by divided difference
operators, are also found within the same algebra $\Nic_W$.

Our point of view of the Nichols\dash Woronowicz algebras is via 
the braided group theory developed by Majid in mid\dash 1990s (see e.g.\
\cite{Mexamples,Mcalculus,Mbook}). 
The principal idea of Majid's approach is that to each object in a
braided category, there is canonically associated a pair of Hopf
algebras in this braided category, which are non-degenerately dually paired. 
We call each of these dually paired braided Hopf algebras a
Nichols\dash Woronowicz (sometimes simply Nichols) algebra.

The term `Nichols algebra' was introduced by Andruskiewitsch and
Schneider in \cite{AS1} and refers to an equivalent definition of this
object as a graded braided Hopf algebra, generated by its degree one
component which is the set of primitives. These conditions first
appeared in the work \cite{N} of Nichols.  
Apparently the first explicit construction of a Nichols\dash
Woronowicz algebra 
per se appeared in the paper  \cite{Wo} by Woronowicz, where  
exterior algebras for quantum differential calculi were studied; 
the term `Woronowicz exterior algebra' is used
by a number of authors. The Nichols\dash Woronowicz algebras seem to
become an increasingly popular object of study.
 
The relations in the Nichols algebra are Woronowicz relations which
ensure the non\dash degeneracy of the duality pairing mentioned above. 
They are not as explicit as the relations in $\mathcal{E}_n$ or
$\mathit{BE}(W,S)$, and in practice, we do not work with the relations
directly. 
We use the methods of braided differential calculus, which makes our
approach essentially different from that of \cite{FK} and \cite{KM}.

How to apply these
methods to the Fomin\dash Kirillov algebras in the case of symmetric
group, was explicitly shown by Majid in \cite{M}. In particular,
the divided difference operators are interpreted as restrictions of 
braided partial derivatives --- the version of this for an arbitrary
Coxeter group plays a central role in Section \ref{fifthsection} below.
It was also proposed in \cite{M} to replace the algebras
$\mathcal{E}_n$ by their Woronowicz quotient and to extend the
construction to other Coxeter groups, which is achieved in the
present paper.

Our construction was inspired, besides the papers mentioned above, 
by the work \cite{MS} of Milinski and Schneider.
In \cite{MS}, a general scheme for constructing Nichols algebras  over a
Coxeter group is 
 discussed (our algebra $\Nic_W$ fits into this scheme), and 
 the Nichols algebra $\Nic_{\Symm_n}$ is explicitly introduced. 
Let us also
 mention a more recent preprint \cite{KMnew} where the `super' 
Nichols\dash Woronowicz  
algebras $\Lambda_w(W)$, which control the noncommutative geometry
 of Weyl groups $W$, are considered. 

The structure of the paper is as follows. In Section
\ref{firstsection} we recall basic facts about Coxeter groups, their
root systems, coinvariant algebras, nilCoxeter algebras and Schubert classes. 
The Nichols algebras are
defined in Section \ref{thirdsection}, after a brief exposition of
braided differential calculi in Section \ref{secondsection}. Our
principal example of the Nichols algebra is $\Nic_W$, a graded braided
Hopf algebra in the Yetter\dash Drinfeld category over 
the Coxeter group $W$, described in
Section~\ref{fourthsection}. 

After all this preparatory work,
in Section~\ref{fifthsection}
we state a principal result of the paper, Theorem \ref{maintheorem},
which implies that $\Nic_W$ contains a graded subalgebra isomorphic to the
coinvariant algebra $S_W$. This commutative subalgebra is generated
by such a subspace $U$ of the degree~$1$ in $\Nic_W$ that
(i)  
$U$ is isomorphic, as a $W$\dash module, to the reflection
representation of $W$; (ii) $U$ is `generic'. 
We describe all such $U$, which we call
generic reflection submodules. 

It turns out that if $W$ is a crystallographic Coxeter group with a simply
laced Dynkin diagram, there is exactly one canonical reflection
submodule, which generates a canonical copy of $S_W$ in $\Nic_W$. If
$W=\Symm_n$, this reflection submodule is precisely the space
generated by Dunkl elements defined in \cite{FK}. Note that we 
do not use explicit expressions for specific Dunkl elements; 
that they generate a copy of $S_W$ follows solely from the fact that
their span is a reflection representation of $W$. For a Coxeter group
of a non\dash simply laced type, there always is more than one way to
embed $S_W$ into $\Nic_W$, but any embedding may be obtained from a
given one by composing with an explicit automorphism of $\Nic_W$.

Section~\ref{nilcoxeter_section} is devoted to the realisation of the
nilCoxeter algebra $\nilcoxeter$ in $\Nic_W$. 
Thus,
two different objects, the coinvariant algebra $S_W$ and the
nilCoxeter algebra $\nilcoxeter$, are identified with subalgebras in the
Nichols\dash Woronowicz algebra $\Nic_W$; 
the (right) action of $\nilcoxeter$ on $S_W$ is interpreted as the
natural action of the braided 
Hopf algebra $\Nic_W$ on itself by derivations; the non\dash
degenerate pairing between $S_W$ and $\nilcoxeter$ is just the
self\dash duality pairing on $\Nic_W$ restricted to these two
subalgebras.

The most important observation in Section~\ref{nilcoxeter_section} is that the
simple root generators in the Nichols\dash Wo\-ro\-no\-wicz algebra $\Nic_W$
obey the Coxeter relations. This fact, interesting by itself, is proved
by expressing the Woronowicz symmetriser of both sides of a Coxeter
relation in terms of paths in the Bruhat graph of a Coxeter group.

In the last section of the paper we mention that $\Nic_W$ is
a quotient of the Fomin-Kirillov algebra $\mathcal{E}_n$ 
when $W=\Symm_n$, and show that $\Nic_W$ is a quotient of the
Kirillov\dash Maeno 
bracket algebra for many other Coxeter groups $W$. We finish by
repeating a conjecture 
made by several authors, that $\Nic_{\Symm_n}=\mathcal{E}_n$.

It should be noted that there are $q$\dash deformed versions of
Fomin\dash Kirillov and Kirillov\dash Maeno algebras (see \cite{FK}
and \cite{KM}), which are intended to be a model for small quantum
cohomology rings of flag manifolds. 
A more recent preprint \cite{KMquantum} by Kirillov and Maeno, which
uses the results of the preliminary version of the present paper,
suggests an embedding of a quantum cohomology ring for a
crystallographic Coxeter group $W$ into a tensor square of the
Nichols\dash Woronowicz
algebra $\Nic_W$.   
We have, however, left the issue of quantising the Nichols\dash
Woronowicz algebra model,  
as well as questions related to the structural
theory of Nichols algebras such as finite dimensionality and Hilbert
series, beyond the scope
the present paper.

\textbf{Acknowledgments.}
I am grateful to Shahn Majid, who introduced me to a broad circle of topics
including braided Hopf algebras and noncommutative geometry. 
This paper owes much to communication with 
Arkady Berenstein, whose remarks and suggestions were very useful.
I thank Victor Ginzburg, Istv\'an Heckenberger and Hans\dash J\"urgen
Schneider for a number of
stimulating discussions during the conference dedicated to the memory of Joseph
Donin, Haifa, July 2004, and Sergey Fomin for his comments on a
preliminary version of this paper. 
I am grateful to 
 Anatol N.\ Kirillov and Toshiaki Maeno for a number of insightful detailed 
discussions on bracket algebras and Nichols\dash Woronowicz algebras.
The work was supported by EPSRC grant GR/S10629.

\section{The coinvariant algebra of $W$}

\label{firstsection}


\subsection{The root system and the reflection representation}
Let $W$ denote a finite Coxeter group 
generated by $s_1,\dots,s_r$ subject to the relations
$(s_i s_j)^{m_{ij}}=1$ $(1\le i,j \le r)$,
where the integers $m_{ij}$ satisfy 
$m_{ii}=1$, $m_{ij}=m_{ji}\ge 2$ for $i\ne j$.
The length of an element $w$ of $W$, denoted by $\len(w)$, 
is defined as the smallest possible number $l$ of factors 
in a decomposition $w=s_{i_1}s_{i_2}\dots s_{i_l}$ of $w$ into a
product of Coxeter generators $s_i$; any such decomposition with
$l=\len(w)$ is called reduced. 

Let $\h$ be a vector space
with a fixed basis $\alpha_1, \dots, \alpha_r$ 
and a non-degenerate symmetric bilinear form $(\cdot, \cdot)$ 
defined by $(\alpha_i, \alpha_j)=-\cos(\pi/m_{ij})$. 
To each $\alpha\in\h$ satisfying $(\alpha, \alpha)=1$ there is
associated an orthogonal reflection 
$h \mapsto h-2(h,\alpha)\alpha$ of $\h$. 
Let the generators $s_i$ of $W$
act on $\h$  
by the reflections associated to $\alpha_i$; this gives rise to the reflection
representation $W\to \mathrm{GL}(\h)$, which is faithful. 
The action of $W$ preserves the bilinear form $(\cdot, \cdot)$ on
$\h$.

The vectors $\alpha_1,\dots, \alpha_r$ are simple roots; 
all $W$\dash images of the $\alpha_i$ in $\h$ are roots and
form the root system $R$ of $W$.
The construction implies that $(\alpha,\alpha)=1$ for all $\alpha\in R$. 
A root which can be written as $\sum_{i=1}^r c_i\alpha_i$, with
nonnegative real $c_i$, is called positive.   
One has $R=R^+\sqcup -R^+$ where $R^+$ is the set  
of positive roots.

If $\alpha$ is a root, that is, $\alpha=w(\alpha_i)$ for some $w\in
W$,
$1\le i \le r$,
then $s_\alpha := w s_i w^{-1}$ acts on $\h$ as the reflection
associated to $\alpha$. Therefore, $s_\alpha=s_{-\alpha}$ is a
well\dash defined element of $W$ (which does not depend on the choice
of $w$ and $i$).     
 
This construction 
of the root system and the reflection (also called
geometric) representation of $W$  
is given in Part II of \cite{Hu}. The space $\h$ can be defined
over the field of real numbers; 
we nevertheless consider its complexification and
assume the ground field to be $\field$. 
The reflection representation of $W$ 
is irreducible, if and only if $W$ is irreducible as a
Coxeter group \cite[V.\S4.7]{Bou}. 
If $W$ is a crystallographic Coxeter group \cite[VI.\S2.5]{Bou},
the reflection
representation $\h$ of $W$ may be identified with a Cartan subalgebra of a
complex semisimple Lie algebra $\mathfrak{g}$, and $W$ with the Weyl
group of $\mathfrak{g}$.  

The $W$\dash action on $\h$ extends to the symmetric algebra
$S(\h)=\bigoplus_{n\ge 0} S^n(\h)$ of $\h$. 
We will refer to the elements of $S(\h)$ as polynomials. 

\subsection{The coinvariant algebra of $W$}
\label{coinv_alg}

By a fundamental result of Chevalley \cite{C},
the subalgebra $S(\h)^W$ of $W$\dash invariant polynomials in $S(\h)$ 
is itself a
free commutative algebra of rank $r$. 
Its generators $f_1,\dots,f_r$ may be chosen to be 
homogeneous polynomials of degrees $\deg f_i=m_i+1$, where 
$1\le m_1\le \dots\le m_r$ are integers
which depend on $W$ but not on the choice of a particular set of generators.
These $m_i$ 
are called the exponents of the Coxeter group $W$. 

Let $I_W$ be the ideal in $S(\h)$ generated by
$f_1,\dots f_r$; in other words, $I_W = S(\h)S(\h)^W_+$ where
$S(\h)^W_+$ is the set of $W$\dash invariant polynomials without
constant term.
The coinvariant algebra of $W$ is, by definition, 
\begin{equation*}
               S_W = S(\h) / I_W.
\end{equation*}
The ideal $I_W$ is $W$\dash stable and graded, 
hence $S_W$ is a graded $W$\dash module.
As shown in \cite{C}, there is an ungraded module isomorphism between
$S_W$ and the regular representation of $W$.
In particular, the dimension of $S_W$ is equal to the number of
elements in $W$.  
See \cite{BL} for information on the graded module structure of $S_W$.


\subsection{The nilCoxeter algebra of $W$}
\label{nilcoxeter_def}
The nilCoxeter algebra $\nilcoxeter$ of $W$ arises naturally in
connection with the coinvariant algebra $S_W$. Namely, $\nilcoxeter$
is the algebra generated by the divided difference operators acting on
polynomials $f\in S(\h)$ as well as on their classes in $S_W$. We
now recall the abstract definition of $\nilcoxeter$ and its
representation by divided difference operators.

Let $\nilcoxeter$ be the algebra generated by $r$ generators $u_1,\dots,u_r$
subject to the nilCoxeter relations
\begin{equation*}
          u_{i}u_{j}u_{i}\ldots 
          = u_{j}u_{i}u_{j}\ldots
          \ \text{($m_{ij}$ factors on each side);} 
\qquad
           u_{i}^2 
           = 0.
\end{equation*}
For $w\in W$ with a reduced decomposition $w=s_{i_1}\dots s_{i_l}$, 
define $u_{w}$ as 
the product $u_{i_1}\dots u_{i_l}$.
(The element $u_{w}$ is well\dash defined because of the
Coxeter relations between the $u_i$.)
The elements $u_{w}$, $w\in W$, are known to form a homogeneous linear basis
of the algebra $\nilcoxeter$; the multiplication table of
$\nilcoxeter$ in terms of this basis is  
\begin{equation*}
u_{v}u_{w}=
\bigg\{
\begin{matrix}
u_{vw}, &\text{ if }\len(v)+\len(w)=\len(vw);
\\
0,\hfill             &\text{ if }\len(v)+\len(w)>\len(vw).
\end{matrix} 
\end{equation*}
Thus, $\nilcoxeter$ is a graded algebra of dimension $|W|$.

\subsection{Divided difference operators}
\label{demazure_operators}

Let $\alpha$ be a root, and $s_\alpha\in W$ be the corresponding reflection.
The linear operator 
\begin{equation*}
    \demazure{\alpha} \colon S(\h) \to S(\h), 
    \qquad
    \demazure{\alpha} f = \frac{f - s_\alpha(f)}{\alpha},
\end{equation*}
is called the divided difference operator. 
(These operators were introduced independently by 
Bernstein\dash Gelfand\dash Gelfand and by  Demazure). 
One may note that the polynomial $f - s_\alpha(f)$ is always divisible by 
$\alpha$, so that the rational function 
$\demazure{\alpha}f$ is a polynomial.
Since the ideal $I_W\subset S(\h)$ is preserved by the action of
$\demazure{\alpha}$, the operators $\demazure{\alpha}$ act on the
coinvariant algebra $S_W$ as well. 


Write $\demazure{i}$ for the divided difference operator 
$\demazure{\alpha_i}$
corresponding to a simple root $\alpha_i$ $(1\le i \le r)$. 
By \cite[IV.\S1--\S2]{H}, the operators $\demazure{i}$ satisfy the
nilCoxeter relations \ref{nilcoxeter_def}; this gives rise to a
representation of the nilCoxeter algebra $\nilcoxeter$ 
on $S(\h)$ and on $S_W$, with the generator $u_i$ acting as
$\demazure{i}$. Either of these representations is faithful.

\subsection{The right action of $\nilcoxeter$ on the coinvariant  algebra}
\label{nilcoxeter_pairing}

We have described, following Hiller \cite{H}, a left action of the
nilCoxeter algebra $\nilcoxeter$ on $S_W$. We will now
convert it to a right action, using the algebra isomorphism 
between $\nilcoxeter$ and its opposite algebra
$\nilcoxeter^\mathrm{op}$ given on the linear basis by 
$u_w \mapsto  u_{w^{-1}}$.
A right action is more convenient for the purpose of
Nichols\dash Woronowicz algebra realisation of
$\nilcoxeter$ and $S_W$.

Define the right\dash hand divided difference operators 
$\rdemazure{\alpha}$ by 
\begin{equation*}
      f\rdemazure{\alpha} = \demazure{\alpha}f, 
 \quad
      f\rdemazure{i} = \demazure{i}f, 
\quad f\in S_W,
\end{equation*}
and let $\nilcoxeter$ act on $S_W$ from the right via 
\begin{equation*}
     f u_i = f \rdemazure{i}.
\end{equation*}
In terms of the basis $\{u_w\}$ of $\nilcoxeter$, this
right action is given by $fu_w = \demazure{w^{-1}}(f)$. 

Observe that the operator $\rdemazure{\alpha}$ obeys the following
version of twisted Leibniz rule: 
\begin{equation*}
          (f g)\rdemazure{\alpha} = 
           f\cdot (g\rdemazure{\alpha}) + (f \rdemazure{\alpha}) s_\alpha(g)
\end{equation*}
for $f$, $g\in S_W$, cf.\ \cite[IV.\S1]{H}. That is,
$\rdemazure{\alpha}$ is a $(1,s_\alpha)$\dash twisted derivation of
the algebra $S_W$.

\subsection{The non-degenerate pairing between $S_W$ and $\nilcoxeter$}
\label{coinv_pairing}

Let $\epsilon \colon S_W\to \field$ be the projection to the degree
zero component $S_W^0=\field$ of $S_W$ (so, $\epsilon f$ is the
`constant term' of $f$). 
Consider a bilinear pairing 
\begin{equation*}
     \langle \cdot,\cdot\rangle_{S_W,\nilcoxeter} \colon
     S_W \otimes \nilcoxeter \to \field, 
     \qquad 
     \langle f, u_w \rangle = \epsilon(f\rdemazure{w}).
\end{equation*}
Similarly to \cite[IV.\S1]{H}, this pairing is non\dash degenerate.
This, in particular, means that for any non\dash
zero $f\in S_W$ there exists a product 
$\rdemazure{} = \rdemazure{i_1}\dots \rdemazure{i_l}$ $(l\ge 0)$
of right\dash hand divided difference operators, such that
$0\ne f \rdemazure{}  \in  \field$ (this will be used later in
\ref{lemma_kernel}).


\subsection{Schubert classes and Schubert polynomials}

In the crystallographic case, the Schubert classes in the cohomology
of the flag variety correspond to a distinguished 
linear basis of $S_W$. This basis $\{\bar X_w \mid w\in W\}$,
described by Bernstein, Gelfand and Gelfand, and
independently by Demazure, 
may be defined purely algebraically as follows: 
$\bar X_{w_0} = \frac{1}{|W|}\prod_{\gamma\in R^+} \gamma$ where $w_0$ is
the longest element in $W$, and for an arbitrary $w\in W$ one has
$\bar X_w = \demazure{w^{-1}w_0} \bar X_{w_0}$.
In fact, $\{\bar X_w \mid w\in W\}$
is the basis of $S_W$ which is dual to $\{u_w\mid
w\in W\}\subset \nilcoxeter$ with respect to the above duality pairing 
$\langle \cdot, \cdot \rangle_{S_W, \nilcoxeter}$. 

The term `Schubert polynomials' usually refers to a family of elements
$X_w\in S(\h)$ which project to $\bar X_w\in S_W$ and satisfy as many
combinatorial properties of $\bar X_w$ as possible.  
For the symmetric group $W=\Symm_n$, the Schubert polynomials were
introduced by Lascoux and Sch\"utzenberger, see the survey in 
\cite{Mac}; a construction of Schubert
polynomials for other
classical types was later suggested independently by
Billey and Haiman in \cite{BH} and by Fomin and Kirillov in \cite{FKschub}. 
Schubert polynomials are intended to be a realisation of $S_W$ in a
free commutative algebra, whereas the approach of \cite{FK}, \cite{KM}
and the present paper is to view $S_W$ as a subalgebra in a
noncommutative algebra. 


\section{Free braided differential calculus}

\label{secondsection}

In this section we recall the free braided differential calculus, 
as  introduced by
Majid e.g.\ in \cite{Mcalculus}. The proofs can be found in Chapters 9
and 10
of \cite{Mbook}.

\subsection{Braided Hopf algebras}
\label{braided_hopf}

Recall that a braiding in a tensor category $(\mathcal{C},\otimes)$ is a 
functorial
family of isomorphisms $\Psi_{A,B}\colon A\otimes B \to B\otimes A$
(where $A$, $B$ are any two objects in $\mathcal{C}$), satisfying the
`hexagon axioms' 
$(\id_B\otimes \Psi_{A,C})(\Psi_{A,B}\otimes
\id_C)${}$=\Psi_{A,B\otimes C}$ and $(\Psi_{A,C}\otimes
\id_B)(\id_A\otimes \Psi_{B,C})${}$=\Psi_{A\otimes B, C}$.
The associativity isomorphisms such as $A\otimes (B\otimes C)\cong
(A\otimes B)\otimes C$, which allow us not to pay attention to the
order of brackets  
in multiple tensor products, are suppressed but are part of the tensor
category setup.
A braided category, as defined in \cite{JS}, is a tensor category
equipped with a braiding.

Let us recall the definition of a braided Hopf algebra,
as given e.g.\ in \cite{Mexamples} (see also a self\dash
contained exposition in \cite{Mbook}). 
Suppose $A$, $B$ are algebras in a braided tensor
category $(\mathcal{C},\otimes,\Psi)$, meaning that their product morphisms 
$\mult_A\in \Hom(A\otimes A,A)$, 
$\mult_B\in \Hom(B\otimes B,B)$
are fixed, as well as the unit morphisms 
$\eta_A\in \Hom(\mathbb{I}, A)$, 
$\eta_B\in \Hom(\mathbb{I}, B)$, where $\mathbb{I}$ is the unit object
in $\mathcal{C}$.

The tensor product of $A$ and $B$ in
$\mathcal{C}$ is also equipped with an algebra structure in the following way:
\begin{equation*}
     \mult_{A\otimes B} = (\mult_A\otimes \mult_B)\circ 
     (\id_A \otimes \Psi_{B,A}\otimes \id_B) 
      \colon A\otimes B \otimes A\otimes B \to  A \otimes B;
\qquad
     \eta_{A\otimes B} = \eta_A \otimes \eta_B.
\end{equation*}
The resulting
braided tensor product algebra  is denoted by
$A\underline\otimes B$.

A braided bialgebra is an object $B$ in $(\mathcal{C},\otimes,\Psi)$ 
which is an algebra and a coalgebra (with coproduct $\Delta$ and 
counit $\epsilon$) such that
$\Delta \colon B \to B\underline\otimes B$ is a morphism of algebras. 
A braided Hopf algebra is a braided bialgebra with 
antipode $S\colon B\to B$. Note that the antipode is braided\dash
antimultiplicative, i.e.\ $S\circ \mult = \mult\circ \Psi_{B,B}\circ
(S\otimes S)$.

There is also a standard notion of graded braided Hopf algebra, meaning 
that $B=\oplus_{n\ge 0} B^n$ in $\mathcal{C}$, 
and the structure morphisms $\mult$, $\eta$, $\Delta$, $\epsilon$, $S$ of $B$ 
respect this grading.

\subsection{Free braided groups}

A braided linear space $(V,\Psi)$ is a pair consisting of a linear
space $V$ and a linear operator $\Psi\colon V\otimes V \to V\otimes
V$ obeying the braid equation
\begin{equation*}
      \Psi_{12}\Psi_{23}\Psi_{12} = \Psi_{23}\Psi_{12}\Psi_{23} \in
      \End V^{\otimes 3}.
\end{equation*}
We use the standard `leg notation' $\Psi_{12}$ etc.\ for the action of
matrices on tensor powers. The braiding $\Psi$ is assumed to come from a braided
category, where the braid equation is a consequence of the hexagon
axioms. 

Suppose $V$ to be finite\dash dimensional and $\Psi$ to be invertible.
Consider the full tensor
algebra $T(V)$ with braiding canonically extended from $V$ by the
hexagon axioms. 
The coproduct $\Delta\colon T(V)\to T(V)\underline\otimes T(V)$,
counit $\epsilon \colon T(V)\to \field$ and antipode 
$S\colon T(V)\to T(V)$ are defined by their values on generators: 
\begin{equation*}
    \Delta v = v\otimes 1 + 1\otimes v, \qquad
    \epsilon v = 0, \qquad
    S v = -v 
\qquad \text{for }v\in V,
\end{equation*}
which makes $T(V)$ a graded braided Hopf algebra, called a `free
braided group'.

To the linear dual $V^*$
of $V$ with the braiding $\Psi^*\in \End ({V^*}^{\otimes 2})
=\End ((V^{\otimes 2})^*)$
there corresponds the graded braided Hopf algebra $T(V^*)$.
We denote the coalgebra maps of $T(V^*)$ by the same letters 
$\Delta$, $\epsilon$, $S$ and use the Sweedler notation 
$\Delta a = a_{(1)}\otimes a_{(2)}$.

\subsection{Braided duality pairing}
\label{pairing}

The evaluation pairing $\langle \xi, v \rangle = \xi(v)$ 
between $V^*$ and $V$ may be extended to a pairing 
\begin{equation*}
  \langle \cdot, \cdot \rangle \colon T(V^*) \otimes T(V) \to \field
\end{equation*}
satisfying the axioms of duality pairing of (braided) Hopf algebras
\cite[(5)]{Mquasi}:
\begin{gather*}
           \langle \phi\psi, x \rangle = \langle \phi, x_{(2)} \rangle 
                                   \langle \psi, x_{(1)} \rangle,
\qquad
           \langle \phi, xy \rangle = \langle \phi_{(2)}, x \rangle 
                                   \langle \phi_{(1)}, y \rangle,
\\
  \langle 1, x \rangle = \epsilon x,
 \qquad \langle \phi, 1 \rangle = \epsilon \phi, 
 \qquad \langle S \phi , x \rangle = \langle \phi, S x \rangle.
\end{gather*}
This braided duality pairing depends on $\Psi$ and does not coincide
with the standard pairing between tensor powers.  
However, one necessarily has $\langle {V^*}^{\otimes m}, V^{\otimes
  n}\rangle=0$ unless $m=n$.

\subsection{Braided partial derivatives}
\label{derivatives}

The space $V^*$ acts on $T(V)$ via left braided partial
derivatives \cite{Mquasi}, \cite[10.4]{Mbook}: 
\begin{equation*}
          \xi \in V^* \quad \mapsto \quad D_\xi \colon T(V)\to T(V),
\qquad 
           D_\xi f = \langle \xi, f_{(1)} \rangle f_{(2)}.
\end{equation*}
Similarly, one has a right action of $V$ on $T(V^*)$ defined by
\begin{equation*}
           v \in V \quad \mapsto \quad \rderiv{v} \colon T(V^*)\to T(V^*),
\qquad 
           \phi \rderiv{v} = \phi_{(1)} \langle \phi_{(2)}, v \rangle .
\end{equation*}
In the present paper, we use right derivatives $\rderiv{v}$ (Remark
\ref{choice_deriv} below explains this chioce and discusses another 
possile choice of derivatives). Here are  
some of the properties of the $\rderiv{v}$.
The operator $\rderiv{v}$ on $T(V^*)$ is adjoint to the multiplication by $v$
from the left in $T(V)$, i.e.\  
 $\langle \phi \rderiv{v}, x\rangle = \langle \phi, vx \rangle$.
For $x=v_1\otimes v_2\otimes \dots \otimes v_m \in
V^{\otimes m}$, denote $\phi \rderiv{x}=\phi \rderiv{v_1}\dots \rderiv{v_m}$.
One can then obviously rewrite the braided duality pairing as
\begin{equation*}
      \langle \phi, x \rangle = \epsilon (\phi \rderiv{x}).
\end{equation*}

\subsection{Braided Leibniz rule}
\label{braided_Leibniz}
One may use an equivalent definition of operators $\rderiv{v}$ via the 
condition $\xi \rderiv{v}=\xi(v)\in T^0(V^*)=\field$ for $\xi\in
T^1(V^*)=V^*$ and the braided Leibniz rule
\begin{equation*}
           (\phi\psi)\rderiv{v} = \phi (\psi\rderiv{v}) + 
           \phi \,\Psi^{-1}_{V,T(V^*)} (\psi \otimes \rderiv{v}),
\qquad \phi,\psi\in T(V^*).
\end{equation*}
Let us explain the notation used in this formula. 
The operator $\Psi_{V,T(V^*)}$ is the braiding between $V$ and
$T(V^*)$ in the braided category $\mathcal{C}$. 
It can be expressed solely in terms of the braiding $\Psi_{V,V}$ on $V$
\cite[Proposition 10.3.6]{Mbook}, however, we will not use this
explicit expression in the general case.
(In \ref{braided_Leibniz_alpha}, we will use the braided Leibniz rule
for a particular braided space $V$ with  
a simple formula for $\Psi_{V,T(V^*)}$.)
Now
  $\phi \, \Psi^{-1}_{V,T(V^*)}(\psi \otimes
\rderiv{v})$ is interpreted in the following way:
one applies the inverse braiding $\Psi^{-1}_{V,T(V^*)}$ 
to $\psi\otimes v$ and obtains an element, say 
$\sum_i v_i \otimes \psi_i$, of $V\otimes T(V^*)$;
one then computes $\sum_i (\phi \rderiv{v_i})\psi_i \in T(V^*)$.

\subsection{Braided symmetriser}
\label{braided_symmetriser}

The above duality pairing $\langle \cdot, \cdot \rangle$ can be written down
explicitly as follows. 

Let $\Braid_n$ denote the braid group with braid generators
$\sigma_1,\dots,\sigma_{n-1}$, and let $\Symm_n$ be the corresponding
symmetric group generated by Coxeter generators $s_1,\dots, s_{n-1}$. 
The Matsumoto section $t\colon \Symm_n \to \Braid_n$ 
is a set\dash theoretical map 
defined by the rule $t(\pi)=\sigma_{i_1}\sigma_{i_2}\dots
\sigma_{i_l}$, whenever $\pi=s_{i_1}s_{i_2}\dots s_{i_l}$ is a reduced
decomposition of $\pi\in \Symm_n$. 
The element $\Sigma_n=\sum_{\pi\in\Symm_n}t(\pi)\in \field\Braid_n$ 
is called the braided (or quantum) symmetriser.

For $k=1,2,\dots$, the braided integers  
$[k]_\sigma\in \field \Braid_n$, the `shifted' braided
integers
$[k]_\sigma^{(s)} \in \field \Braid_n$, 
and the braided factorials 
$[k]!_\sigma\in \field\Braid_n$ are:
\begin{align*}
&[k]_\sigma 
=1+\sigma_1 + \sigma_2\sigma_1+\dots +\sigma_{k-1}\dots \sigma_2\sigma_1,
\\
&[k]_\sigma^{(s)} 
=1+\sigma_{s} + \sigma_{s+1}\sigma_{s}+\dots +\sigma_{s+k-2}\dots \sigma_{s+1}\sigma_{s},
\\
&[k]!_\sigma = 
[k]_\sigma [k-1]_\sigma^{(2)} [k-2]_\sigma^{(3)} \dots [2]_\sigma^{(k-1)},  
\end{align*}
cf.\ \cite{Mcalculus}.
%
The braided symmetriser factorises as $\Sigma_n = [n]!_\sigma$.

Let the generator $\sigma_i$ of $\Braid_n$ act on 
$V^{\otimes n}$ as $\Psi_{i,i+1}$. 
Denote the resulting action of 
braided integers, resp.\ braided factorials,
by $[k]_\Psi$, resp.\ $[k]!_\Psi \in \End V^{\otimes n}$.
Then the duality pairing \ref{pairing} between ${V^*}^{\otimes n}$ and 
$V^{\otimes n}$ is explicitly given by  
\begin{equation*}
          \langle \phi , x \rangle = (\phi \mid [n]!_\Psi x)
            =([n]!_{\Psi^*}\phi \mid  x),           
\end{equation*}
where $(\cdot | \cdot)$ is the evaluation pairing
$(\xi_n\otimes\dots \otimes \xi_2\otimes \xi_1 \mid v_1\otimes
v_2\otimes \dots \otimes v_n)=\prod_{i=1}^n \xi_i(v_i)$.

\section{Nichols-Woronowicz algebras}

\label{thirdsection}

\subsection{Definition of the Nichols-Woronowicz algebra}

Let $V$ be a linear space with a braiding $\Psi$, and let $T(V)$,
$T(V^*)$ be the braided Hopf algebras introduced in the previous
section. The duality pairing $\langle \cdot, \cdot \rangle \colon
T(V^*)\times T(V) \to \field$ may be degenerate.
Let $I(V^*)$, resp.\ $I(V)$, be the kernel  of the pairing in
$T(V^*)$, resp.\ $T(V)$. 
Since we deal with a duality pairing between graded (braided) Hopf algebras,
the kernels $I(V^*)$,  $I(V)$ are graded Hopf ideals.

The algebras
\begin{equation*}
      \Nic(V^*) = T(V^*)/I(V^*), \qquad \Nic(V) = T(V)/I(V)
\end{equation*} 
are called the Nichols\dash Woronowicz algebras of $V^*$ and $V$.

We now state the basic properties of the Nichols\dash Woronowicz
(also called Nichols) algebra which follow
directly from the construction outlined above. 
For a braided
space $V$, this construction 
 leads to a dual pair of Nichols algebras $\Nic(V^*)$,
$\Nic(V)$; we formulate the properties for $\Nic(V^*)$, keeping in
mind that their analogues hold for $\Nic(V)$ as well.

\begin{xlemma}
\label{basic_properties}
$(i)$ $\Nic(V^*)=\oplus_{n\ge 0}\Nic(V^*)^n$ is a graded braided Hopf algebra.

$(ii)$ $\Nic(V^*)^0=\field$, $\Nic(V^*)^1=V^*$.

$(iii)$ There is a Hopf algebra duality pairing
  $\langle\cdot, 
  \cdot \rangle\colon \Nic(V^*)\otimes \Nic(V) \to \field$, which is
  non\dash degenerate.

$(iv)$ $\Nic(V^*)$ is generated by $\Nic(V^*)^1$ as an algebra.
\end{xlemma}
\begin{proof}
The ideal $I(V^*)$ is graded, hence $(i)$; the duality pairing
$\langle\cdot, \cdot \rangle$ between $T(V^*)$ and $T(V)$
 is non\dash degenerate in $T^0(V^*)=\field$
and $T^1(V^*)=V^*$, hence $(ii)$. The meaning of the construction of
$\Nic(V^*)$ is that one eliminates the kernel of the duality pairing,
so $(iii)$ follows; 
$(iv)$ is obvious since $\Nic(V^*)$ is a quotient of $T(V^*)$.
\end{proof}

\subsection{Woronowicz relations}
\label{woronowicz}

Description \ref{braided_symmetriser} of the duality pairing $\langle
\cdot, \cdot \rangle$ means that, in degree $n$, the kernel of the pairing is
precisely the kernel of the braided symmetriser $\Sigma_n$.
Thus we arrive at the following presentation of the Nichols algebras:
\begin{equation*}
  \Nic(V^*) = \bigoplus_{n\ge 0} {V^*}^{\otimes n}/\ker [n]!_{\Psi^*},
\qquad
  \Nic(V) = \bigoplus_{n\ge 0} V^{\otimes n}/\ker [n]!_{\Psi}.
\end{equation*}
This presentation (with $-\Psi$ instead of $\Psi$ which does not
affect the braid equation) was used as the definition of the exterior algebra
for quantum differential calculi in the work of Woronowicz \cite{Wo}.

The relations in $\Nic(V^*)$, $\Nic(V)$, although given as kernels of
quantum (anti) symmetrisers,  are in general not known explicitly and
may be complicated. 
One may try 
to find the rank of $[n]!_\Psi$ which is the dimension of $\Nic(V)^n$,
or to check if $\Nic(V)$ is at all finite\dash dimensional,
but this is usually difficult 
even when the braided space $V$ is `small'.
An excellent example of this kind of work
in the case $\dim V=2$ is the paper \cite{He} of
Heckenberger.

\subsection{An equivalent definition of the Nichols algebra}

The earliest occurrence of the following properties, 
characterising the Nichols algebra, was in
the work of Nichols \cite{N}:
(1) $\Nic(V)^0=\field$, $\Nic(V)^1=V$;
(2) $\Nic(V)^1=P(\Nic(V))$;
(3) $\Nic(V)^1$ generates $\Nic(V)$ as an algebra.
Here $P(\Nic(V))$ is the set of primitive elements of $\Nic(V)$ 
(i.e.\ those $a\in \Nic(V)$ 
satisfying $\Delta a =a\otimes 1 + 1\otimes a$).
A proof of equivalence between this definition and, say, Woronowicz
presentation  \ref{woronowicz} can be found in \cite{Sch}.

\subsection{Braided partial derivatives}
\label{b_p_deriv}

For $v\in V$, one has the braided partial derivative  
$\rderiv{v}\in \End T(V^*)$ 
which satisfies $\langle \phi\rderiv{v}, x\rangle = \langle \phi, vx
\rangle$. It follows that if $\phi\in T(V^*)$ lies in the kernel of
the pairing $\langle\cdot,\cdot \rangle$, then $\phi\rderiv{v}$ does.
Therefore, $\rderiv{v}$ are well\dash defined endomorphisms of 
the Nichols\dash Woronowicz algebra $\Nic(V^*)$. 

The right action of vectors $v\in V$ on $\Nic(V^*)$ by braided partial
derivatives $\rderiv{v}$ extends to
a right action of the algebra $\Nic(V)$ on $\Nic(V^*)$.
In fact, the operators $\rderiv{x}$ defined for arbitrary
$x\in\Nic(V)$ 
as $\phi \rderiv{x} = \phi_{(1)}\langle \phi_{(2)}, x \rangle$, satisfy
\begin{equation*}
        \phi \rderiv{xy} = (\phi \rderiv{x})\rderiv{y}.  
\end{equation*}

The  duality pairing $\langle \cdot, \cdot \rangle$
between $\Nic(V^*)$ and $\Nic(V)$ is induced from $T(V^*)$, $T(V)$,
and is therefore given by 
$\langle \phi, x \rangle = \epsilon(\phi \rderiv{x})$, where 
$\phi\in\Nic(V^*)$, $x\in \Nic(V)$.
Non\dash degeneracy of this pairing obviously implies that 
the right action of $\Nic(V)$ on $\Nic(V^*)$ is faithful.  

The following criterion, which describes the joint kernel
of all braided partial derivatives in $\Nic(V^*)$, turns out to be
extremely useful 
when working with Nichols\dash Woronowicz algebras. 
It asserts that a
`function' , all of whose partial derivatives are zero, is a constant. 
The criterion is equivalent to the non-degeneracy of 
$\langle \cdot, \cdot \rangle$, thus is automatic in the braided differential
calculus given by $\Nic(V^*)$, $\Nic(V)$; 
a form of it in the classical Hopf algebra theory approach
can be traced back to Nichols \cite{N}.

\begin{xcriterion}
\label{criterion_orig}
The following are equivalent: 
(a) $\phi\in\Nic(V^*)$ is a constant.
(b) $\phi \rderiv{v} = 0$ for all $v\in V$.
\end{xcriterion}
\begin{proof}
If $\phi \in \Nic(V^*)^0=\field$, 
then obviously $\phi \rderiv{v}=0$ since $\rderiv{v}$ lowers the degree
by one. Suppose now that a non\dash zero $\phi\in \Nic(V^*)^n$ is in
the kernel of all $\rderiv{v}$. By the non\dash degeneracy of the pairing,
there exists $x\in \Nic(V)^n$ such that $\langle \phi, x\rangle \ne 0$;
since $\Nic(V)$ is generated by $V$ as an algebra, one may choose $x$
to be $v_1v_2\dots v_n$ for some $v_i\in V$.
Then $\epsilon(\phi \rderiv{v_1}\dots \rderiv{v_n})\ne 0$.
But since $\phi\rderiv{v}=0$ for any $v$, this can only be possible if
$n=0$.
\end{proof}

%
%
%
%
%
%
%
%
%
%

\subsection{Two simplest examples}

The two simplest examples of Nichols algebras, for any linear space $V$, are: 
(a) $\mathrm{Sym} (V)$, the symmetric algebra of $V$ 
(let $\Psi=\tau$ to be the flip $\tau(x\otimes y) = y\otimes x$);
(b) $\bigwedge V$, the exterior algebra of $V$ (let $\Psi=-\tau$). 
Both cases are observed easily by \ref{woronowicz}
if one notes that the Woronowicz symmetriser $[n]!_\Psi$ 
becomes the  usual (a) symmetrisation, (b) antisymmetrisation map
on $V^{\otimes n}$.
In the case of symmetric algebra, both left and right braided
derivatives $D_v$, $\rderiv{v}$ coincide with the usual directional
derivative $\frac{\partial}{\partial v}$ on the polynomial ring;
on $\bigwedge^n V$, the left braided derivative $D_v$ is the
contraction (or interior product) operator corresponding to $v$ 
and differs from
$\rderiv{v}$ by a factor of $(-1)^{n-1}$.


\section{Nichols algebras $\Nic_W$ over Coxeter groups}
\label{fourthsection}

In this section, our main example of Nichols\dash Woronowicz algebra is
introduced. The braided category it comes from, is the Yetter\dash
Drinfeld module category over the Coxeter group $W$, whose definition
we now recall.
 
\subsection{The Yetter-Drinfeld module category over a finite group}
Let $\Gamma$ be a finite group. The objects of the 
 Yetter\dash Drinfeld module category $\YD{\Gamma}$ over $\Gamma$
are linear spaces $V$ 
with the following structure:

(1) the $\Gamma$\dash action $\Gamma\times V\to V$, $(g,v)\mapsto gv$;

(2) the $\field\Gamma$\dash coaction, which is the same as
    $\Gamma$\dash grading $V=\oplus_{g\in \Gamma}V_g$;

(3) the compatibility condition $gV_h = V_{ghg^{-1}}$.

One knows that $\YD{\Gamma}$ is a braided tensor category: for
$U,V\in\mathit{Ob}(\YD{\Gamma})$, the $\Gamma$\dash action on
$U\otimes V$ is $g(u\otimes v)=gu\otimes gv$ and the $\Gamma$\dash
grading is $(U\otimes V)_g = \oplus_{h\in \Gamma} U_h\otimes
V_{h^{-1}g}$.
For $V\in\mathit{Ob}(\YD{\Gamma})$, the braiding is given by
$\Psi(x\otimes y) = gy \otimes x$ whenever $x\in V_g$, $y\in V$. 

Of course, the general theory of Nichols algebras applies well for this
particular type of braided linear spaces. 
The Nichols algebras in the Yetter\dash Drinfeld category over a group
have been an
object of extensive study, in particular because they are linked with
pointed Hopf algebras --- see for example the survey \cite{AS} of
Andruskiewitsch and Schneider.

\subsection{The Yetter-Drinfeld module $\VV{W}$}
\label{V_W}

We now specify $\Gamma$ to be the Coxeter group $W$ and will introduce
a particular braided space $\VV{W}\in\mathit{Ob}(\YD{W})$, 
thus linking the
content of Section~\ref{firstsection} with the Nichols algebra theory 
of Sections \ref{secondsection} and \ref{thirdsection}.
We will freely use the notation from all preceding sections.

Let $\VV{W}$ be the linear space spanned by symbols 
$[\alpha]$ where $\alpha$ is a root of $W$, subject to the relation 
$[-\alpha]=-[\alpha]$.   
The dimension of $\VV{W}$ is thus $|R^+|$.

The $W$\dash action on $\VV{W}$ is given by 
$w[\alpha]=[w\alpha]$, and the $W$\dash grading  
is given by assigning the degree $s_\alpha$ to the basis element $[\alpha]$.
The action and the grading are compatible, so
that $\VV{W}$ is a Yetter\dash Drinfeld module over $W$.
The resulting braiding $\Psi$ on $\VV{W}$ is given explicitly 
by 
\begin{equation*}
\Psi([\alpha]\otimes[\beta]) = [s_\alpha \beta] \otimes 
[\alpha].
\end{equation*}
Our main object is the Nichols algebra 
$\Nic(\VV{W})$.

\begin{xremark}
The definition of $V_W$ in fact comes from at least two sources. 

First, precisely this linear space is the degree $1$ component in the
bracket algebra $BE(W,S)$, defined by Kirillov and Maeno in
\cite{KM}. 
The algebra $BE(W,S)$ has the same quadratic relations as
$\Nic(\VV{W})$, but in general, bracket algebras of \cite{KM} have less
relations then $\Nic(\VV{W})$ and are not Nichols algebras.
(See \ref{bracket_algebras} below for more details.)

Second, there is a recipe by Milinski and Schneider
\cite[section 5]{MS}, which, applied to the Coxeter group $W$,
suggests to take the vector space $V$ with basis $\{x_t\}_{t\in T}$ 
where $T$ is the set of all reflections in $W$. The basis element $x_t$ is of
$W$\dash degree $t$.   
The $W$\dash action on $V$ is given by
$gx_t=\chi(g,t)x_{gtg^{-1}}$, where the function 
$\chi\colon W\times T \to \field$
satisfies $\chi(gh,t)=\chi(g,hth^{-1})\chi(h,t)$, so that $V$ is a
Yetter\dash Drinfeld module.

Example 5.3 in \cite{MS} defines $\chi$ in the case $W=\Symm_n$ by
$\chi(g,t) = 1$ if $g(i)<g(j)$, $\chi(g,t)=-1$ if $g(i)>g(j)$, where 
$t=(ij)$ is a reflection, $1\le i<j\le n$, and $g\in \Symm_n$. 
Our Yetter\dash Drinfeld module $V_W$ generalises this example to the
case of an arbitrary Coxeter group: indeed, 
put $x_{s_\alpha}=[\alpha]$ for $\alpha\in R^+$, and define the
function $\chi$ by 
\begin{equation*}
\alpha\in R^+, g\in W \qquad \Rightarrow\qquad       
\chi ( g, x_{s_\alpha} ) = 
\Bigg\{
\begin{matrix}
1, &\text{ if }g(\alpha)\in R^+;
\\
-1,\hfill             &\text{ if }-g(\alpha)\in R^+.
\end{matrix} 
\end{equation*}
\end{xremark}

\subsection{The self-dual Nichols algebra $\Nic_W$}
\label{self-duality}

In the general picture of Section~\ref{thirdsection}, we had two
dually paired Nichols algebras $\Nic(V^*)$ and $\Nic(V)$.
In the case $V=V_W$, however, it is useful to identify $V_W$ with its
dual $V_W^*$ via the non\dash degenerate bilinear form on $V_W$
defined by 
\begin{equation*}
             \langle [\alpha], [\beta] \rangle =
             \delta_{\alpha,\beta}, 
\qquad
              \alpha, \beta\in R^+.
\end{equation*} 
This bilinear form is $W$\dash invariant and is compatible, in a
proper sense, with the
$W$\dash grading on $\VV{W}$. Thus, $\VV{W}$ and $\VV{W}^*$ are
isomorphic as objects in the category $\YD{W}$.
The braiding $\Psi\in\End V_W^{\otimes 2}$ becomes self\dash adjoint
with respect to the evaluation pairing on $V_W^{\otimes 2}$.
The Nichols algebras $\Nic(V_W^*)$ and $\Nic(V_W)$ are then
canonically identified:
\begin{equation*}
    \Nic_W := \Nic(V_W^*) = \Nic(V_W).
\end{equation*}
Note that $\Nic_W$, being a braided Hopf algebra in $\YD{W}$, is a
$W$\dash module. The $W$\dash action is given by
$w([\alpha_1][\alpha_2]\dots [\alpha_n]) = [w\alpha_1][w\alpha_2]\dots
[w\alpha_n]$. 

%

\subsection{Braided partial derivatives in $\Nic_W$}
\label{bderiv_alphabeta}

To each positive root $\alpha$ there corresponds, by the general
construction outlined above and the self\dash duality of $\Nic_W$, the
right braided partial derivative $\rderiv{[\alpha]}$ on $\Nic_W$. 
Let us restate the main properties of braided partial derivatives for
this particular Nichols algebra. 

First, the restriction of $\rderiv{[\alpha]}$ onto $V_W=\Nic_W^1$ is
given by
\begin{equation*}
        [\beta] \rderiv{[\alpha]}=
\Bigg\{
\begin{matrix}
\pm 1, &\alpha = \pm\beta,
\\
0,  &\alpha \ne \pm \beta,
\end{matrix} 
\qquad\beta\in R.
\end{equation*}

\subsection{}
\label{braided_Leibniz_alpha}

Second, braided Leibniz rule \ref{braided_Leibniz} simplifies:
since $V_W$ is self\dash dual as was shown \ref{self-duality}, 
the braiding $\Psi_{V_W,V_W^*}$ is now the same as
$\Psi_{V_W,V_W}\colon [\alpha]\otimes f \mapsto s_\alpha(f)\otimes
[\alpha]$ for $\alpha\in R$, $f\in V_W$. 
The inverse braiding is thus given by 
$\Psi^{-1}_{V_W,V_W}(f\otimes [\alpha]) = [\alpha]\otimes
s_\alpha(f)$. 
It follows from the hexagon
axiom that $\Psi^{-1}_{V_W, T(V_W)}$ is given by the same
formula (but with $f\in T(V_W)$); passing to the quotient, one obtains
the same inverse braiding between $V_W$ and $\Nic_W$. The braided Leibniz rule
now becomes
\begin{equation*}
      (fg)\rderiv{[\alpha]} =
      f(g\rderiv{[\alpha]}) + (f\rderiv{[\alpha]}) s_\alpha(g), 
\qquad f,g\in \Nic_W, \alpha\in R.
\end{equation*}
That is, $\rderiv{[\alpha]}$ is an $(1,s_\alpha)$\dash twisted
derivation of $\Nic_W$.

Finally, \ref{criterion_orig} now gives the following

\begin{xcriterion}
\label{criterion}
$f\in \Nic_W$ is a constant, if and only if $f\rderiv{[\alpha]}=0$ for
all $\alpha\in R^+$.
\end{xcriterion}


\section{The realisation of the coinvariant algebra in~$\Nic_W$} 
\label{fifthsection}

In this section we describe a graded subalgebra
of $\Nic_W$ isomorphic to the coinvariant algebra
of the Coxeter group $W$. 
Thus, the Nichols algebra $\Nic_W$ provides
a model for the coinvariant algebra, Schubert calculus 
and (in the crystallographic case) cohomology of the flag manifold
for an arbitrary Coxeter group $W$ in the same sense as the Fomin\dash
Kirillov algebras $\mathcal{E}_n$ from \cite{FK} provide such a
model for $W=\Symm_n$. The relationship between $\Nic_W$ and the
 Fomin\dash Kirillov algebras will be discussed in the next Section.

The degree\dash preserving $W$\dash equivariant embedding 
$S_W \hookrightarrow \Nic_W$ turns out
to be unique up to a composition with certain automorphisms of $\Nic_W$,
which we describe explicitly.

\subsection{Reflection submodules in $\VV{W}$}

The algebra $S_W$ is generated by its degree $1$ component
$S^1(\h)=\h$, which is the reflection representation of $W$.
Hence, a graded subalgebra of $\Nic_W$ which is isomorphic to $S_W$ as
a graded algebra, must be generated by $U\subset \VV{W}=\Nic_W^1$,
such that $U$ is a $W$\dash 
submodule of $\VV{W}$ isomorphic to the reflection representation $\h$
of $W$.  

We will, however, be slightly more general and consider all 
non\dash zero submodules
$U\subset \VV{W}$ which are images of $W$\dash homomorphisms 
$\mu \colon \h \to \VV{W}$. 
Such submodules $U$ of $\VV{W}$ will be called reflection
submodules.

\subsection{The support of a submodule}
\label{support}

Before we describe subalgebras generated by reflection submodules,
let us introduce a bit more notation. 
Let $U$ be a linear subspace of $\VV{W}$. 
Define the support of $U$ by
\begin{equation*}
        \supp U = \{ \alpha \in R \mid U \rderiv{[\alpha]}\ne 0\}.
\end{equation*}
In other words, the support of $U$ is the minimal set of
roots $\pm\alpha$, such that the linear span of $[\alpha]$ in $\VV{W}$
contains $U$.   
If $U$ is a $W$\dash invariant subspace of $\VV{W}$, 
then $\supp U$ is a $W$\dash invariant subset of $R$, and therefore is
itself a root system in $\h$. Let $W(\supp U)\subseteq W$ be  
the group generated by
reflections with respect to the roots in $\supp U$.

\subsection{Generic reflection submodules}
\label{generic_def}
We call a submodule $U\subset\VV{W}$ generic, if 
$\supp U$ is the whole of $R$. 
The generic reflection submodules are singled out by the condition
$W(\supp U)=W$.  We justify the term `generic' later in \ref{generic}.

The rest of this section will mainly be devoted to the proof 
of the following 
\begin{xtheorem}
\label{maintheorem}
$(i)$ Let $U$ be a reflection submodule in $\VV{W}$.  
The subalgebra generated by $U$ in $\Nic_W$ 
is commutative.
It is isomorphic,
as a graded algebra, to the coinvariant algebra of the Coxeter group
$W(\supp U)$.

$(ii)$
Generic reflection submodules of $V_W$ exist.
Each such submodule 
generates a subalgebra of $\Nic_W$ isomorphic to $S_W$.
\end{xtheorem}

Let us start with a Lemma which provides an explicit description of
reflection submodules in $\VV{W}$. 
\begin{xlemma}
\label{formula_mu}
$(1)$ Any $W$\dash module homomorphism 
$\mu \in \Hom_W(\h, \VV{W})$ is given by a formula 
\begin{equation*}
         \mu(x) = \sum_{\alpha\in R} c_\alpha(x,\alpha)[\alpha],
\end{equation*}
where $\alpha \mapsto c_\alpha$ is a $W$\dash invariant scalar
function on the root system $R$. 

$(2)$ The support of $\mu(\h)$ is $\{\alpha\in R \mid c_\alpha\ne 0\}$.

$(3)$ $\mu$ is injective if and only if $\supp \mu(\h)$ spans $\h$.
\end{xlemma}
\begin{proof}
$(1)$ 
Consider a new $W$\dash module $\widetilde{\VV{W}}$, with linear basis 
$\{v_{\alpha} \mid \alpha\in R\}$ and the $W$\dash action given by 
$wv_\alpha = v_{w\alpha}$. 
The module $\VV{W}$ is a submodule of $\widetilde{\VV{W}}$, via the inclusion 
$[\alpha] = v_\alpha-v_{-\alpha}$. 
For any linear map $\mu$ from $\h$ to $\widetilde{\VV{W}}$, there are 
elements $b^\alpha$ of $\h$ 
such that $\mu(x) = \sum_{\alpha\in R} (x,b^\alpha)v_\alpha$.
The map $\mu$ is $W$\dash equivariant, if and only if 
$b^{w\alpha}=wb^\alpha$ for any root $\alpha$ and any element $w$ of
the group  
$W$.

Now $\Hom_W(\h,\VV{W})$ will consist of those $W$\dash maps $\mu\colon \h \to
\widetilde{\VV{W}}$ whose image lies in $\VV{W}$. In terms of the elements
$b^\alpha$ this translates to $b^{-\alpha} = -b^\alpha$. 
On the other hand, let $s_\alpha\in W$ be the reflection associated to
the root $\alpha$; since $s_\alpha \alpha = -\alpha$, 
one has $s_\alpha b^\alpha = b^{-\alpha} = -b^\alpha$. 
This immediately implies that $b^\alpha$ is proportional to
$\alpha$, say $b^\alpha = c_\alpha \alpha$,
and the $W$\dash equivariance of the $b^\alpha$ implies that 
$c_{w\alpha}=c_\alpha$ for any element $w$ of $W$;
%
%
so $(1)$ follows. $(2)$ is immediate from the definition \ref{support} of
support. The kernel of $\mu$ consists of those $x\in\h$ which are 
orthogonal to all $c_\alpha \alpha$, i.e.\ $\ker \mu = (\supp\mu(\h))^{\perp}$, hence $(3)$.  
\end{proof}


\subsection{}
\label{generic}
We now have the following information on reflection 
submodules in $\VV{W}$, immediate from Lemma \ref{formula_mu}. 
Let, say, $\mathit{RS}$ be the variety of all reflection submodules of
$\VV{W}$, and let $\mathit{RS}_{\cong
  \h}$ (resp.\ $\mathit{RS}_{\text{\it gen}}$) be the part of 
$\mathit{RS}$ consisting of submodules isomorphic to $\h$
(resp.\ generic submodules).
First of all, $\mathit{RS}$ is not empty and 
is of dimension equal to the number of
$W$\dash orbits in $R$ minus one. 
Furthermore, 
$\mathit{RS}_{\cong \h} \supseteq \mathit{RS}_{\text{\it gen}}$;
the generic part 
$\mathit{RS}_{\text{\it gen}}$, as well as $\mathit{RS}_{\cong \h}$, 
is an open dense set in $\mathit{RS}$ (any $W$\dash invariant 
function $\alpha\mapsto c_\alpha$ on $R$, such that $c_\alpha\ne 0$ for all $\alpha$, gives rise to a generic reflection submodule). In particular,
generic reflection
submodules of $\VV{W}$ exist.

\subsection{The multiplicity of $\h$ in $\VV{W}$}
If $W$ is an irreducible Coxeter group \cite[IV.\S1.9]{Bou}, 
the reflection representation $\h$ 
is irreducible \cite[V.\S4.7-8]{Bou}.
The multiplicity of $\h$ in $\VV{W}$ is
then equal to the number of $W$\dash orbits in the root system $R$.
If, moreover, $W$ is a Weyl group of simply laced type (so that there
is only one orbit in $R$), then
there is a canonical non\dash zero reflection submodule in $\VV{W}$,
which is generic.

If $W=\Symm_n$, this canonical reflection submodule is precisely the
subspace spanned by Dunkl elements in the terminology of \cite{FK}.
 
Our next step is to establish the commutativity of subalgebras in
$\Nic_W$ generated by reflection submodules. 

\begin{xproposition}
\label{proposition_commute} 
Let $U$ be a reflection submodule of $\VV{W}$.
The subalgebra $\langle U \rangle$ of $\Nic(\VV{W})$, generated by
$U$, is commutative.
\end{xproposition}
\begin{proof}
Let $U$ be the image of a $W$\dash module 
map $\mu \colon \h \to \VV{W}$.
We will show that any two elements of $U$ commute in $\Nic_W$.
By the formula for $\mu$ given in Lemma \ref{formula_mu}, 
two elements of $U$ can be written as 
$\mu(x)=\sum_{\alpha\in R}
c_\alpha(x,\alpha^\vee)[\alpha]$ and similarly $\mu(y)$.
The commutator $[\mu(x), \mu(y)]$ is an element of degree $2$ in
$\Nic_W$.
According to presentation \ref{woronowicz} of $\Nic_W$, 
the commutator vanishes if and only if 
\begin{equation*}
(\id + \Psi)(\mu(x)\otimes \mu(y)) = (\id+\Psi)(\mu(y)\otimes
\mu(x)),
\end{equation*}
where $\Psi$ is the braiding \ref{V_W} of $\VV{W}$.
The left hand side rewrites as 
\begin{equation*}
\sum_{\alpha,\beta\in R} c_\alpha c_\beta 
\bigl((x,\alpha)(y,\beta)+(x,\beta)(y,s_\beta
\alpha) \bigr)
[\alpha]\otimes[\beta];
\end{equation*}
since $s_\beta \alpha = \alpha - 2(\alpha,
  \beta)\beta$, 
this equals 
\begin{equation*}
\sum_{\alpha,\beta\in R} c_\alpha c_\beta 
\bigl((x,\alpha)(y,\beta)+(x,\beta)(y,
\alpha) 
-2(x,\beta)(y,\beta)(\alpha,\beta)
\bigr)
[\alpha]\otimes[\beta].
\end{equation*}
This expression for the left hand side is symmetric in $x$ and $y$, 
therefore is equal to the right hand side. 
\end{proof}

It follows from the last Proposition that any $W$\dash homomorphism $\mu
\colon \h \to \VV{W}$ extends to a map $\mu \colon S(\h) \to \Nic_W$
of $W$\dash module algebras. 
The kernel of $\mu$ will be calculated using the vanishing
criterion \ref{criterion}, but for that we need to know how to apply
the braided derivations $\rderiv{[\alpha]}$ to $\mu(f)$, $f\in
S(\h)$. 
The following key
Lemma, which is ideologically the same as 
Proposition 9.5 from \cite{FK},
shows how to express $\mu(f)\rderiv{[\alpha]}$
in terms of the divided difference operator 
$\rdemazure{\alpha}$ acting on $S(\h)$.

\begin{xlemma}
\label{lemma_Dd}
Suppose that an algebra homomorphism 
$\mu\colon S(\h) \to \Nic_W$ is defined by $\mu(x) =
\sum_{\beta\in R}c_\beta(x,\beta)[\beta]$ for $x\in\h$. 
Then
\begin{equation*}
\mu(f) \rderiv{[\alpha]} = c_\alpha \mu(f \rdemazure{\alpha})
\qquad\text{for }f \in S(\h), \ \alpha\in R.
\end{equation*}
\end{xlemma}
\begin{proof}
The maps $F_1(f)=\mu(f)\rderiv{[\alpha]}$,
$F_2(f)=c_\alpha\mu(f
\rdemazure{\alpha})$ from $S(\h)$ to $\Nic_W$ 
vanish on constants.
Apply them to $x\in\h=S^1(\h)$.
By \ref{bderiv_alphabeta},
$F_1(x)=c_\alpha(x,\alpha)[\alpha]+c_{-\alpha}(x,-\alpha)[-\alpha]
=2c_\alpha(x,\alpha)$;
since $x\rdemazure{\alpha}=2(x,\alpha)$, 
one has $F_2(x)=2c_\alpha(x,\alpha)$,
hence $F_1$ and $F_2$ agree on $\h$.
By \ref{braided_Leibniz_alpha} and 
\ref{demazure_operators},
both are extended to products of elements of $\h$
according to the twisted Leibniz rule 
$F_i(fg)=\mu(f)F_i(g)+F_i(f)\mu(s_\alpha(g))$.
Therefore, $F_1=F_2$.
\end{proof}
\begin{xremark}
\label{choice_deriv}
It is this lemma that explains why we chose the right partial
derivatives $\rderiv{[\alpha]}$ over their seemingly more convenient
left\dash hand counterparts $D_{[\alpha]}$. The advantage of
$\rderiv{[\alpha]}$ is 
that this satisfies the $(1,s_\alpha)$\dash twisted Leibniz rule, and
as a consequence, coincides, up to a scalar factor, with the divided
difference operator $\demazure{\alpha}$ on $S_W$. The left derivatives
$D_{[\alpha]}$ on $\Nic_W$ do not obey such a reasonable Leibniz
rule. 

Another choice of braided partial derivatives to realise the
divided difference operators would be the braided\dash left
derivatives $\bar D_\xi$, $\xi\in V^*$ which are defined on $T(V)$,
for an arbitrary braided space $V$, by $\bar D_\xi f 
= (\id_{T(V)}\otimes \langle\cdot, \cdot \rangle)(\Psi_{V^*,T(V)}(\xi \otimes f_{(1)})\otimes f_{(2)})$. 
The derivatives $\bar D_{[\alpha]}$ were used
in the case $W=\Symm_n$ in \cite{M}. The operator $\bar D_{[\alpha]}$ 
on $\Nic_W$ satisfies the
$(s_\alpha,1)$\dash twisted Leibniz rule which is equally good for the
divided difference operator $\demazure{\alpha}$. The only drawback for
us is that the $\bar D_{[\alpha]}$ do not give rise to a representation of 
$\Nic_W$ on itself. They lead to an action of a new algebra $\tilde \Nic_W$,
with the same underlying linear space as $\Nic_W$ but with twisted
multiplication $f\star g = \cdot\circ \Psi^{-1}_{\Nic_W, \Nic_W}(f\otimes
g)$. 
Using $\bar D_{[\alpha]}$ instead of $\rderiv{[\alpha]}$, one can
modify the content of this and the next Sections, replacing 
$\Nic_W$ with its twisted version $\tilde\Nic_W$ where necessary.
\end{xremark}
\begin{xcorollary}
\label{corollary_vanish}
Let $\mu\colon S(\h) \to \Nic_W$ be as above, and let 
$W'=W(\supp\mu(\h)) \subseteq W$. Then $\mu(f)=0$ for
 any homogeneous $W'$\dash invariant polynomial $f\in S(\h)$ of
 positive degree.
\end{xcorollary}
\begin{proof}
Take a homogeneous $f\in S(\h)^{W'}_+$.  
If a root $\alpha$ is in $\supp \mu(\h)$ so that $s_\alpha$ is in $W'$,
then $s_\alpha(f)=f$ and 
$f\rdemazure{\alpha}=0$, therefore $\mu(f)\rderiv{[\alpha]}=0$ by Lemma
\ref{lemma_Dd}. If $\alpha\not\in\supp\mu(\h)$, then 
$c_\alpha=0$ by Lemma \ref{formula_mu}(2), hence 
$\mu(f)\rderiv{[\alpha]}=0$ again by Lemma
\ref{lemma_Dd}. Thus, $\mu(f)$ lies in the kernel of all
$\rderiv{[\alpha]}$,
which implies that $\mu(f)\in \field$ by Criterion \ref{criterion}.
But since $\mu(f)$ is of positive degree, this means that $\mu(f)=0$.
\end{proof}
\begin{xremark}
\label{remark_vanish}
It follows from the Corollary that the kernel of $\mu$ contains the ideal
$I_{W'}(\h):=$ $S(\h)S(\h)^{W'}_+$ of $S(\h)$. Let $U=\mu(\h)$; then 
$\mu$ induces a surjective map 
from $S(\h)/I_{W'}(\h)$ onto the 
subalgebra $\langle U \rangle$ of $\Nic_W$. 
%
%
%
%

Note that $\h$ may not be the reflection representation for $W'=W(\supp U)$ 
because 
$\h'=\mathrm{span}(\supp U)$ is not necessarily the whole of $\h$ 
(although it is, if $W$ is
an irreducible Coxeter group). Still,  
$S(\h)/I_{W'}(\h)$ is
isomorphic to the coinvariant algebra $S_{W'}=S(\h')/I_{W'}$. 
Indeed, $\h=\h'\oplus \mathfrak{k}$ where 
the action of $W'$ on $\mathfrak{k}$ is trivial; 
$S(\h)=S(\h')\oplus J$ and $I_{W'}(\h) = I_{W'} \oplus J$ 
where $J=S(\h')S(\mathfrak{k})_+$,  so that
the isomorphism follows. 

Thus, we have already proved that there is an onto map $S_{W'}\to
\langle U \rangle$. To complete the steps needed for the proof of
Theorem~\ref{maintheorem}, 
we have to show that this map is an isomorphism.
We are going to use lemma \ref{lemma_Dd} one more time.
\end{xremark}
\begin{xlemma}
\label{lemma_kernel}
In the above notation,
the kernel of $\mu\colon S(\h) \to \Nic_W$ is precisely 
$I_{W'}(\h)=S(\h)S(\h)^{W'}_+$.
\end{xlemma}
\begin{proof}
The inclusion $I_{W'}(\h)\subseteq \ker \mu$ has been demonstrated in
the last Corollary and Remark. We assume now that $f\in S(\h)$ does
not lie in $I_{W'}(\h)$ and show that $\mu(f)\ne 0$. Decompose
$S(\h)=S(\h') \oplus J$ as in the Remark; since $J\subseteq
I_{W'}(\h)\subseteq \ker \mu$, it is enough to assume that 
$f\in S(\h')$ (and $f\not\in I_{W'}$). By \ref{coinv_pairing},
there exist roots $\gamma_1,\gamma_2,\dots,\gamma_l$ in the root system $\supp
U$, such that $f \rdemazure{\gamma_1}\rdemazure{\gamma_2}\dots
\rdemazure{\gamma_l} \in a+I_{W'}$  
where $a\in S^0(\h)$ is a 
non\dash zero constant. Then by Lemma \ref{lemma_Dd} one has 
$\mu(f)\rderiv{[\gamma_1]}\rderiv{[\gamma_2]}\dots
\rderiv{[\gamma_l]} = c_{\gamma_1}c_{\gamma_2}\dots
c_{\gamma_l}a$, which is not zero since $c_{\gamma_i}\ne 0$
by Lemma \ref{formula_mu}(2). Hence $\mu(f)\ne 0$. 
\end{proof}
\subsection{Proof of the Theorem}
The proof of Theorem \ref{maintheorem} is already contained in 
\ref{formula_mu}--\ref{lemma_kernel}, but we summarise it here for clarity.
$(i)$ 
If $U\subset\VV{W}$ is a reflection submodule, write $U=\mu(\h)$
where $\mu\colon \h \to \VV{W}$ is a $W$\dash module map, and extend
this map by \ref{proposition_commute} to a surjective homomorphism $\mu \colon
S(\h) \to \langle U \rangle$ of algebras. By Lemma \ref{lemma_kernel},
the kernel of $\mu$ in $S(\h)$ is $I_{W'}(\h)$ where $W'=W(\supp
U)$. Therefore, $\langle U \rangle$ is isomorphic to 
$S(\h)/I_{W'}(\h)$, which is 
$S_{W'}$ as shown in \ref{remark_vanish}. 
$(ii)$ Generic reflection submodules exist by \ref{generic}. 
If $U$ is such a submodule, i.e.\ $\supp U=R$ and $W(\supp U)=W$, 
then  $\langle U \rangle \cong S_W$ by part $(i)$. 
Theorem \ref{maintheorem} is proved.

\subsection{Automorphisms of the Nichols\dash Woronowicz algebra
which permute copies of $S_W$ in $\Nic_W$}
\label{permute}

We have seen that, if there is more than one $W$\dash orbit in the
root system $R$ of $W$, a degree\dash preserving embedding of $S_W$
into $\Nic_W$ is not unique. 
Such an embedding is determined by assigning a non\dash zero value 
of $c_\alpha$
to each $W$\dash
orbit in $R$.
However, two such
embeddings always differ by an action of an automorphism of
$\Nic_W$. This easily follows from the explicit description of generic
reflection submodules in $\VV{W}$ given in Lemma~\ref{formula_mu}:

\begin{xlemma}
\label{permute_lemma}
Let $\mu, \mu'\colon S_W \hookrightarrow \Nic_W$ be two degree\dash
preserving  embeddings. There is a Hopf algebra automorphism
$\theta\colon \Nic_W \to \Nic_W$ such that $\mu' = \theta\circ\mu$. 
\end{xlemma}
\begin{proof}
We may assume $\mu$ to be fixed; let the restriction of $\mu$ to
$\h=S^1_W$ be given by $\mu(x) = \sum_{\alpha\in
  R}(x,\alpha)[\alpha]$. For an arbitrary embedding $\mu'$ one has 
$\mu'(x) = \sum_{\alpha\in R}c_\alpha (x,\alpha)[\alpha]$ with some
non\dash zero coefficients $c_\alpha$, $W$\dash invariant in $\alpha$.
Define an invertible linear map $\theta\colon \VV{W}\to \VV{W}$ by
$\theta([\alpha])=c_\alpha [\alpha]$. Then $\theta$ is an automorphism
of $\VV{W}$ as a Yetter\dash Drinfeld module over $W$, because $\theta$
is compatible with both the $W$\dash action 
and the $W$\dash
grading on $\VV{W}$. Thus, $\theta$ preserves the braiding on $\VV{W}$
and therefore extends to a Hopf algebra automorphism $\theta\colon
\Nic_W \to \Nic_W$. One has  $\mu'(x) = \theta(\mu(x))$ for $x\in
S^1_W$ and (since both sides are algebra homomorphisms) for all $x\in S_W$. 
\end{proof}

\section{The nilCoxeter subalgebra of $\Nic_W$}
\label{nilcoxeter_section}

In the previous section, we realised the coinvariant algebra $S_W$ of
the Coxeter group $W$ as a subalgebra in the Nichols\dash Woronowicz
algebra $\Nic_W$.

We now recall 
the nilCoxeter
algebra $\nilcoxeter$, which acts on $S_W$ and is
non\dash degenerately paired with $S_W$, as
described in \ref{nilcoxeter_def}--\ref{nilcoxeter_pairing}.
We will now show that all this structure (the nilCoxeter algebra, its
action on $S_W$ and its pairing with $S_W$), and not only the algebra
$S_W$ itself, is realised in $\Nic_W$.

\subsection{}

The subalgebra
of $\Nic_W$ isomorphic to the coinvariant algebra $S_W$, constructed
in Section~\ref{fifthsection}, 
depends on a choice of
a $W$\dash invariant scalar function $\alpha\mapsto c_\alpha\ne 0$ on the
root system $R$. 
From now on, we assume that
\begin{equation*}
           c_\alpha = 1 \quad \text{for all roots $\alpha$,}
\end{equation*}
to simplify the exposition. All results in this section may, however, be
restated for arbitrary $c_\alpha$ by applying the automorphism
$\theta\colon \Nic_W \to \Nic_W$, defined in the proof of Lemma
\ref{permute_lemma}. 

\subsection{}

In light of the above assumption, we consider a linear map
\begin{equation*}
\mu\colon \h \to \VV{W},
\quad
\mu(x) = \sum_{\alpha\in R}(x,\alpha)[\alpha],
\end{equation*}
that extends, as we
know from Theorem~\ref{maintheorem} and its proof, to an
embedding
\begin{equation*}
  \mu\colon S_W \hookrightarrow \Nic_W
\end{equation*}
of algebras. 

By Lemma~\ref{lemma_Dd}, the partial derivative operator $\rderiv{[\alpha]}$ 
corresponding to a root $\alpha$ acts on the subalgebra $\mu(S_W)$ as
the divided difference operator $\rdemazure{\alpha}$. 
Therefore, the restrictions $\rderiv{[\alpha_i]}|_{\mu(S_W)} \in \End
\mu(S_W)$
of the braided partial derivatives
corresponding to
simple roots $\alpha_1,\dots,${}$\alpha_r$,
onto the subalgebra $\mu(S_W)$
satisfy the nilCoxeter relations
\ref{nilcoxeter_def} just as the divided difference operators $\rdemazure{i}$
do. 
It turns out that 
$\rderiv{[\alpha_i]}$ themselves (and not only their restrictions to
the finite\dash dimensional subalgebra $\mu(S_W)$ of $\Nic_W$) satisfy
the nilCoxeter relations. This is shown in the
next Theorem.


Recall from \ref{b_p_deriv} that the operators $\rderiv{x}$ on
the self\dash dual Nichols\dash Woronowicz algebra $\Nic_W$ are
defined  for arbitrary $x\in \Nic_W$. 
If $x=[\alpha][\beta]\dots [\gamma]$ in $\Nic_W$,
$\alpha,\beta,\dots,${}$\gamma\in R$, one has 
$f\rderiv{x}=f\rderiv{[\alpha]}\rderiv{[\beta]}\dots\rderiv{[\gamma]}$ 
for $f\in \Nic_W$.

\begin{xtheorem}
\label{nilcoxeter_theorem}
(i) The simple root generators $[\alpha_1],\dots,${}$[\alpha_r]$ in
  $\Nic_W$ obey the nilCoxeter relations.

(ii) The map $\nu\colon  \nilcoxeter \to \Nic_W$, given on generators by 
  $\nu(u_i)= [\alpha_i]$ and extended multiplicatively to
  $\nilcoxeter$, is an algebra isomorphism between the
  nilCoxeter algebra $\nilcoxeter$ and the subalgebra of
  $\Nic_W$ generated by $[\alpha_1],\dots${}$[\alpha_r]$.

(iii)  The right action of $\nilcoxeter$ on $S_W$ is expressed in terms of
  right derivations of $\Nic_W$: $\mu(f u) =
  \mu(f)\rderiv{\nu(u)}$ for 
  $f\in S_W$, $u \in \nilcoxeter$. 

(iv) The non\dash degenerate pairing $\langle \cdot, \cdot
  \rangle_{S_W,\nilcoxeter}$ coincides
  with the restriction of the self\dash duality pairing $\langle
  \cdot, \cdot \rangle =\langle
  \cdot, \cdot \rangle_{\Nic_W}$ on $\Nic_W$ to the subalgebras
  $\mu(S_W)$ and $\nu(\nilcoxeter)$ :
\begin{equation*}
           \langle f, u 
  \rangle_{S_W, \nilcoxeter}  = \langle \mu(f), \nu(u) 
  \rangle_{\Nic_W}.    
\end{equation*}  
\end{xtheorem}

\begin{xremark}
As an illustration to this result, one may consider two
diagrams: 
\begin{equation*}
\begin{matrix}
      S_W & \otimes & N_W & \to & S_W
\\
      \cap &     & \cap & & \cap
\\
     \Nic_W & \otimes & \Nic_W & \to & \Nic_W
\end{matrix}
\quad\text{and}\quad
\begin{matrix}
      S_W & \otimes & N_W & \to & \field
\\
      \cap &     & \cap & & \Vert
\\
     \Nic_W & \otimes & \Nic_W & \to & \field
\end{matrix}\ ,
\end{equation*}
where the horizontal arrows denote right action resp.\ pairing, and the
inclusions $\cap$ stand for the embeddings $\mu$,  $\nu$ of $S_W$ and
$N_W$ into $\Nic_W$.
The statement of the Theorem means that both diagrams are commutative.
\end{xremark}

\subsection{}
\label{coxeter_braided}
We start the proof of Theorem~\ref{nilcoxeter_theorem} by verifying
the Coxeter relation between the simple root generators $[\alpha_s]$
and $[\alpha_t]$
of $\Nic_W$. In fact, this is the longest part of the proof; 
for all crystallographic root systems, this can be achieved by an
explicit calculation since it is enough to check the cases $m_{st}=2,3,4,6$.
Our argument, however, is valid for any Coxeter group. 
Because of Woronowicz relations in $\Nic_W$,  the
Coxeter relation
$     [\alpha_s] [\alpha_t] [\alpha_s] \ldots${}$ 
     = [\alpha_t] [\alpha_s] [\alpha_t] \ldots$ 
($m_{st}$ factors on each side)
is equivalent to 
\begin{equation*} 
        [m_{st}]!_\Psi ([\alpha_s] \otimes [\alpha_t] \otimes
	[\alpha_s] \otimes \dots)
      = [m_{st}]!_\Psi ([\alpha_t] \otimes [\alpha_s] \otimes
	[\alpha_t] \otimes \dots), 
\end{equation*}
where $[m_{st}]!_\Psi$ is the braided symmetriser defined
in~\ref{braided_symmetriser}. 

To prove this relation, we will express both sides explicitly in terms
of paths in the Bruhat graph of a dihedral group.

\subsection{The dihedral group}
\label{dihedral}

To check relation~\ref{coxeter_braided}, we first note that we may
restrict ourselves to the root subsystem of rank $2$, generated by the
simple roots $\alpha_s$ and $\alpha_t$. This root subsystem will be of
type $I_2(m)$ for some $m\ge 2$, and will consist of positive roots
$\gamma_0=\alpha_s$, $\gamma_1, \dots$, $\gamma_{m-2}$,
$\gamma_{m-1}=\alpha_t$ and negative 
roots $\gamma_{m+i}=-\gamma_i$, $0\le i\le m-1$. 
These roots may be viewed as vectors $\gamma_i = (\cos \frac{i}{m}\pi,
\sin  \frac{i}{m}\pi)$ in the coordinate plane.
Applying the reflection $s_{\gamma_i}$ to a root $\gamma_j$, one
obtains $\gamma_{m+2i-j}$ (the indices are understood modulo $2m$).
The Coxeter group of type $I_2(m)$ is the dihedral group $\mathbb{D}_m$,
which consists of the following elements:
\begin{align*}
      v_0 = \id; 
      \quad v_{l} &= s_{\gamma_0}s_{\gamma_{m-1}}s_{\gamma_0}\ldots 
      \ (l \text{ factors}),\quad 1\le l \le m;
      \\
      v_{-l} &= s_{\gamma_{m-1}}s_{\gamma_0}s_{\gamma_{m-1}}\ldots 
      \ (l \text{ factors}),
      \quad 1\le l \le m;
      \qquad v_{m}=v_{-m}.
\end{align*}
The Coxeter generators of $\mathbb{D}_m$ are $v_1=s_{\gamma_0}$ and 
$v_{-1}=s_{\gamma_{m-1}}$. 
For any $l=0$,$1,\dots$,$m$,
the length of $v_{\pm l}$ in the group $\mathbb{D}_m$ is $l$.
The group $\mathbb{D}_m$ may be viewed as the parabolic subgroup of
the Coxeter group $W$
generated by the simple reflections with respect to $\alpha_s$ and $\alpha_t$;
the length function $\len(\cdot)$ on $\mathbb{D}_m$ is then induced from $W$.

\subsection{The Bruhat graph of the dihedral group}

Recall that the Bruhat graph of a Coxeter group $W$  
is a labelled directed graph, with $W$ as the set of vertices; 
an edge from $w$ to $w'$ exists if and only if $w'=s_\gamma w$ for a
positive root $\gamma$ of $W$, and $\len(w')=\len(w)+1$; this edge is
labelled by the root $\gamma$ and denoted by $w'\xleftarrow{\gamma} w$. 

We need an explicit description of the Bruhat graph of
$\mathbb{D}_m$. This graph may be drawn as a regular $2m$\dash gon
with vertices $v_0$, $v_1,\dots$, $v_m$, $v_{-(m-1)},\dots$, $v_{-1}$
(in this cyclic order) and sides parallel to roots $\gamma_i$. 
The edges of the graph are 
$v_{l+1} \xleftarrow{\gamma_l} v_l$, 
$v_{-(l+1)} \xleftarrow{\gamma_{m-l}} v_{-l}$ for $0\le l \le m-1$
(the sides of the $2m$\dash gon) and 
$v_{l+1} \xleftarrow{\gamma_0} v_{-l}$, 
$v_{-(l+1)} \xleftarrow{\gamma_{m-1}} v_l$ for 
$1\le l \le m-2$ (the diagonals of the $2m$\dash gon parallel to
$\gamma_0$ or $\gamma_{m-1}$). In the Bruhat graph drawn this way,
each edge is labelled by the positive root parallel to this edge.
 
We will consider paths in the Bruhat graph starting at
the vertex $v_0=\id$. A path $\omega$ of length $l$, consisting of edges 
$v_{\pm l} \xleftarrow{\gamma_{i_l}} v_{\pm(l-1)} \dots 
\xleftarrow{\gamma_{i_2}} v_1
\xleftarrow{\gamma_{i_1}} v_0$, will be denoted
by $\omega=(\gamma_{i_l},\dots,\gamma_{i_2}, \gamma_{i_1})$. 

\subsection{The tensor representation of a Bruhat path}
Consider an injective set\dash theoretical map   
\begin{align*}
            \{\text{paths in the Bruhat graph}\} \quad &\xrightarrow{t}  \quad T(\VV{W})
\\
             \omega = (\gamma_{i_l},\dots, \gamma_{i_2}, \gamma_{i_1}) 
              \quad &\mapsto \quad
              t(\omega)=[\gamma_{i_l}]\otimes \dots \otimes
	      [\gamma_{i_2}]\otimes [\gamma_{i_1}].      
\end{align*}
It is convenient to refer to $t(\omega)$ as the tensor representation
of the path $\omega$.

Let $\omega$ be a Bruhat path from the vertex $v_0=\id$ to a vertex $v\in
\mathbb{D}_m$. 
The $W$\dash degree of
$t(\omega)=[\gamma_{i_{\len(v)}}]\otimes \dots \otimes
[\gamma_{i_1}]$,
i.e., the product $s_{\gamma_{i_{\len(v)}}}\dots s_{\gamma_{i_1}}$, is
equal to $v$, the final vertex of $\omega$.
Note that the braiding $\Psi$ is compatible with the $W$\dash grading
on $T(\VV{W})$: $\Psi_{i,i+1}$ leaves the $W$\dash degree intact.
Thus, if we apply the braiding $\Psi$ at positions $i$, $i+1$ to
$t(\omega)$, we will get an element of $T(\VV{W})$ 
which is either a tensor representation of another path from $v_0$ to
$v$ or not a tensor representation of any Bruhat path at all.
For example, $[\gamma_1]\otimes [\gamma_0]$ is the tensor representation
of a Bruhat path from $v_0=\id$ to $v_2=s_{\gamma_0}s_{\gamma_{m-1}}$,
but one has $\Psi([\gamma_1]\otimes [\gamma_0])=-[\gamma_2]\otimes
[\gamma_1]$, which 
obviously does not
correspond to any Bruhat path (when $m\ge 3$) because of the minus sign. 

\subsection{$\Psi$-generating paths}

The path $\omega$ from $v_0=\id$ to $v\in \mathbb{D}_n$ will be called
a $\Psi$\dash generating path, if the 
Woronowicz symmetrisation of the tensor representation of $\omega$ is
the sum of tensor representations of all Bruhat paths from $v_0$ to $v$:
\begin{equation*} 
       [\len(v)]!_\Psi t(\omega) 
          = \sum_{\omega' = v\leftarrow \dots \leftarrow v_0} t(\omega')       
\end{equation*}
Note that there are $2^{\len(v)-1}$ Bruhat paths from $v_0$ to $v$,
since there are two choices for each intermediate vertex 
$v_{\pm 1},\dots,v_{\pm\len(v)-1}$ of the path. Thus, there are
$2^{\len(v)-1}$ terms on the right hand side of this equation. 
A priori, the left hand side has $\len(v)!$ terms; hence, equality for
$\len(v)>2$ may be possible only due to cancellations on the left.

Although the definition of a $\Psi$\dash generating Bruhat path makes
sense for any Coxeter group $W$, 
in general we can only conjecture that $\Psi$\dash
generating paths exist for any vertex $v$ of the Bruhat graph of $W$.
However, the case of the dihedral group $\mathbb{D}_m$ is handled more
easily because of a very explicit description of the Bruhat graph. We
have the following

\begin{xlemma}
\label{lemma_generating}
(a)
The path $\omega_l^+ = (\gamma_0$, $\gamma_{m-1}$, $\gamma_0$, $\gamma_{m-1},\dots)$
of length $l$ is a $\Psi$\dash generating path from the vertex
$v_0=\id$ to the vertex $v_l$
in the Bruhat graph of the dihedral group $\mathbb{D}_m$.
(b)
The path $\omega_l^- = (\gamma_{m-1}$, $\gamma_0$, $\gamma_{m-1}$, $\gamma_0,\dots)$
of length $l$ is a $\Psi$\dash generating path from $v_0$ to $v_{-l}$. 
\end{xlemma}
\begin{proof}
Denote by $P_l$ the sum of tensor representations 
of all Bruhat paths from $v_0$ to $v_l$. 
Throughout this proof, we are going to write $\gamma_i$ instead of 
$[\gamma_i]$ for the basis elements of $\VV{W}$; this does not lead to a
confusion but ensures that the generators of $T(\VV{W})$ do not mix up with
the braided integers.
 
We have to show that $[l]!_\Psi t(\omega_l^\pm) = P_{\pm l}$.
Induction in $l$; the case $l=0$ is trivial.
If $l=1$, $\omega_1^+=(\gamma_0)$ is the only
  Bruhat path from 
  $v_0=\id$ to $v_1=s_{\gamma_0}$, and, trivially, it is $\Psi$\dash
  generating. Similarly for $\omega_1^-=(\gamma_{m-1})$.

Assume that $l\ge 2$ and that the Lemma is proved for $l-1$ and $l-2$.
The properties 
$[l]!_\Psi = [l]_\Psi(\id\otimes [l-1]!_\Psi)$ and  
$[l]_\Psi=\id+(\id\otimes [l-1]_\Psi)\Psi_{12}$
of braided integers and braided factorials follow from their
definition~\ref{braided_symmetriser}.
One therefore has 
\begin{align*}
     [l]!_\Psi t(\omega_l^+) 
&= [l]!_\Psi \bigl(\gam{0}\otimes
     t(\omega_{l-1}^-) \bigr)
= [l]_\Psi \bigr( \gam{0}\otimes [l-1]!_\Psi t(\omega_{l-1}^-) \bigr)
\\ & =[l]_\Psi \bigl(\gam{0}\otimes P_{-(l-1)}\bigr)
\\ & = \gam{0}\otimes P_{-(l-1)} + (\id\otimes [l-1]_\Psi)\Psi_{12}
\bigl(\gam{0}\otimes P_{-(l-1)} \bigr). 
\end{align*}
Any Bruhat path $\omega$ from $v_0$ to $v_{-(l-1)}$  passes either
through 
the vertex $v_{-(l-2)}$ or through $v_{l-2}$, and the last edge of
$\omega$ is labelled by $\gamma_{m-l+1}$ or by $\gamma_{m-1}$,
respectively. Therefore, $P_{-(l-1)} = \gam{{m-l+1}}\otimes
P_{-(l-2)}
+ \gam{{m-1}} \otimes P_{l-2}$.
Using this, we rewrite
\begin{multline*}
[l]!_\Psi t(\omega_l^+) = \gam{0}\otimes P_{-(l-1)}
\\
+
(\id\otimes [l-1]_\Psi) \bigl(\Psi(\gam{0}\otimes
\gam{{m-l+1}})\otimes P_{-(l-2)} 
+\Psi(\gam{0}\otimes \gam{{m-1}})\otimes P_{l-2}
\bigr).
\end{multline*}
We compute $\Psi(\gam{0}\otimes
\gam{{m-l+1}}) = \gam{{l-1}}\otimes \gam{0}$ 
and $\Psi(\gam{0}\otimes \gam{{m-1}}) = \gam{1}\otimes
\gam{0}$. The tensors $P_{\pm (l-2)}$ are replaced, by the induction
hypothesis, with $[l-2]!_\Psi t(\omega_{l-2}^\pm)$. 
Thus we obtain
\begin{align*}
     [l]!_\Psi t(\omega_l^+)& = \gam{0}\otimes P_{-(l-1)}
+ \gam{{l-1}}\otimes [l-1]_\Psi \bigl( \gam{0}\otimes 
          [l-2]!_\Psi t(\omega_{l-2}^-) \bigr)
\\
      & \qquad + \gam{1} \otimes  [l-1]_\Psi \bigl( \gam{0}\otimes
          [l-2]!_\Psi t(\omega_{l-2}^+) \bigr)
\\ 
& = \gam{0}\otimes P_{-(l-1)}
+ \gam{l-1}\otimes [l-1]!_\Psi \bigl( \gam{0}\otimes t(\omega_{l-2}^-) \bigr)
+ \gam{1}\otimes [l-1]!_\Psi \bigl( \gam{0}\otimes t(\omega_{l-2}^+) \bigr).
\end{align*}
Note that  $\gam{0}\otimes t(\omega_{l-2}^-)=t(\omega_{l-1}^+)$,
so that by the induction hypothesis,
the second term is equal to $\gam{l-1}\otimes P_{l-1}$.
The tensor $\gam{0}\otimes t(\omega_{l-2}^+)$ in the third term is of
the form $\gam{0}\otimes\gam{0}\otimes\ldots$, and lies in the kernel
of the Woronowicz symmetriser $[l-1]!_\Psi$ (indeed,
$[\gamma_0]\cdot [\gamma_0]\cdot \ldots$ is zero in the Nichols\dash
Woronowicz algebra $\Nic_W$). Therefore, the third term on the right
hand side vanishes, yielding $[l]!_\Psi t(\omega_l^+) = \gam{0}\otimes
P_{-(l-1)}+ \gam{l-1}\otimes P_{l-1}$.  
But this is equal to $P_l$, because
a path from $v_0$ to $v_l$ in the Bruhat graph of $\mathbb{D}_m$
either passes via $v_{-(l-1)}$ and has the last edge labelled by
$\gamma_0$, or passes via $v_{l-1}$ and has the last edge labelled by
$\gamma_{l-1}$. 

An argument establishing the other equality $[l]!_\Psi t(\omega_l^-) =
P_{-l}$, is completely analogous. The Lemma is thus proved.
\end{proof}


\subsection{Proof of the Coxeter relations}
We are now ready to prove the Coxeter relation for the simple root
generators $[\alpha_s]=[\gamma_0]$ and $[\alpha_t]=[\gamma_{m-1}]$ in
the Nichols\dash Woronowicz algebra $\Nic_W$. 
Let us show that relation~\ref{coxeter_braided}, which we
rewrite as 
\begin{equation*}
         [m]!_\Psi([\gamma_0] \otimes [\gamma_{m-1}]\otimes [\gamma_0]
         \otimes \dots)
= 
         [m]!_\Psi([\gamma_{m-1}] \otimes [\gamma_0]\otimes [\gamma_{m-1}]
         \otimes \dots),
\end{equation*}
holds. Indeed, since $v_m=v_{-m}$ (both are equal to the longest word in the
group $\mathbb{D}_m$), Lemma~\ref{lemma_generating} implies that 
both sides are equal  to
the sum of tensor representations of all paths from $v_0$ to $v_m$ in
the Bruhat graph of $\mathbb{D}_m$. The Coxeter relation is proved.

\subsection{The rest of the proof of Theorem \ref{nilcoxeter_theorem}} 
To establish part (i) of the
Theorem, it now remains to add that $[\alpha_i][\alpha_i]=0$ in $\Nic_W$
because $(\id+\Psi)([\alpha_i]\otimes [\alpha_i])=0$. 

By (i), there is a well\dash defined algebra homomorphism $\nu\colon
\nilcoxeter \to \Nic_W$ defined by $\nu(u_{i})=[\alpha_i]$.
Let $u=u_{i_1}u_{i_2}\dots u_{i_l}$ be a basis
element
of $\nilcoxeter$.
Then $\rderiv{\nu(u)}$ is
$\rderiv{[\alpha_{i_1}]}\dots \rderiv{[\alpha_{i_l}]}$, and by
Lemma~\ref{lemma_Dd}
$\mu(f)\rderiv{\nu(u)}=\mu(f u)$, so part (iii) of
the theorem follows. 
Part (iv) also follows because $\langle f, u
\rangle_{S_W,\nilcoxeter}$ is the constant term of $fu$ in
$S_W$, which is equal to the constant term of $\mu(fu)$ in
$\Nic_W$; the latter is $\epsilon(f\rderiv{\nu(u)})$ which equals
$\langle \mu(f), \nu(u)\rangle_{\Nic_W}$ by \ref{b_p_deriv}.
We are left to prove part (ii); but (iv) implies that the image of
$\nu$ is non\dash degenerately paired with the $|W|$\dash dimensional
subalgebra $\mu(S_W)$ in $\Nic_W$, therefore $\dim(\im \nu)=|W|$ and
$\nu$ is one\dash to\dash one. 
This finishes the proof of Theorem~\ref{nilcoxeter_theorem}.

\begin{xremark}
Although we proved Coxeter relations between generators of the
Nichols\dash Wo\-ro\-no\-wicz algebra $\Nic_W$ corresponding to simple
roots, the same method shows that any two generators $[\alpha]$ and
$[\beta]$ of $\Nic_W$, $\alpha,\beta\in R$, obey a Coxeter relation up
to a sign.
Indeed, consider the root subsystem of type
$I_2(m)$ generated by $\alpha$ and $\beta$, 
where $m\ge 2$ is such that the scalar product $(\alpha,\beta)$
equals $-c \cos \frac{\pi}{m}$, $c=\pm 1$.
The roots $\alpha, c\beta$ may be chosen as the simple roots in this
root subsystem; the above argument allows one to compute the
Woronowicz symmetriser of $[\alpha]\otimes [c\beta] \otimes
[\alpha]\otimes \ldots$ and 
yields the Coxeter relation of degree $m$ between $[\alpha]$ and
$c[\beta]$.

In particular, if $\gamma$ is the highest root of a crystallographic
root system $R$, the generators
$[\alpha_1],\dots,[\alpha_r],[-\gamma]$ obey the affine nilCoxeter
relations and generate a subalgebra in $\Nic_W$ which is a quotient of the
(infinite\dash dimensional) `affine nilCoxeter algebra'. However, this
quotient is proper and finite\dash dimensional in known cases. This
observation is due to A.~N.~Kirillov. 
\end{xremark}


\section{The algebras $\Nic_W$ and the constructions of
Fomin-Kirillov   and Ki\-ril\-lov-Ma\-eno}

We conclude the paper by outlining the relationship between the
Nichols algebra $\Nic_W$ which we constructed for an arbitrary Coxeter
group $W$, the quadratic algebra $\mathcal{E}_n$ constructed in
\cite{FK} for the symmetric group $\Symm_n$, and the generalisation
$\mathit{BE}(W,S)$ of $\mathcal{E}_n$ for an arbitrary Coxeter group,
defined in \cite{KM}.

\subsection{The quadratic algebra $\Nic_\qd(V_W)$}

Let $\Psi\colon V_W\otimes V_W \to V_W\otimes V_W$ be the braiding on
the Yetter\dash Drinfeld module $V_W$ defined in
Section~\ref{fourthsection}, and let $T(V_W)$ be the free braided
group. Denote by $I_{\qd}(V_W)$ the two\dash sided 
ideal of $T(V_W)$ generated by
$\ker(\id+\Psi)\subset V_W^{\otimes 2}$. Put
\begin{equation*}
      \Nic_{\qd}(V_W) = T(V_W) / I_{\qd}(V_W);
\end{equation*}
that is, to define $\Nic_{\qd}(V_W)$, one imposes only the quadratic
Woronowicz relations on $T(V_W)$. The algebra $\Nic_\qd$ is a braided
Hopf algebra in the category $\YD{W}$ with a self\dash duality pairing which may be
degenerate; the Nichols algebra
$\Nic_W$ is a (possibly proper) quotient of $\Nic_\qd(V_W)$.

\subsection{$\Nic_\qd(V_{\Symm_n})$ is the Fomin-Kirillov algebra}

The algebra $\Nic_\qd(V_{\Symm_n})$ is the same as 
the quadratic algebra $\mathcal{E}_n$, introduced by Fomin and
Kirillov in \cite{FK}. This
was independently observed in \cite{MS} and in \cite{M}. 
These algebras coincide as braided Hopf algebras in the Yetter\dash
Drinfeld module category over $\Symm_n$.
For $1\le a<b \le n$, let $ab$ denote the root
$\alpha_a+\alpha_{a+1}+\dots +\alpha_{b-1}$ in the root system of $\Symm_n$;
the operators $\Delta_{ab}\colon \mathcal{E}_n \to \mathcal{E}_n$, 
defined in \cite[Section 9]{FK}, can be viewed as 
braided
partial derivatives 
on 
$\Nic_\qd(V_{\Symm_n})$, as was noticed in \cite{M}. 
It is also shown in \cite{M} that 
the Hopf algebra structure on the `twisted group algebra'
$\mathcal{E}_n\otimes \field \Symm_n$, introduced and studied in
\cite{FP}, can be obtained by Majid's biproduct bosonisation of
$\mathcal{E}_n$.

\subsection{$\Nic_W$ and Kirillov-Maeno bracket algebras}
\label{bracket_algebras}

For a Coxeter group $W$, the bracket algebra $\mathit{BE}(W,S)$, where
$S$ stands for the set of Coxeter generators, is defined in \cite{KM}
as the quotient of the tensor algebra $T(V_W)$ of the linear space
$V_W$ by the following relations (we use our notation from Sections
\ref{firstsection}--\ref{fourthsection}):
 
(1) $[\gamma]^2=0$ for all $\gamma\in R$;

(2) For any intersection $R'$ of  a $2$\dash dimensional plane in $\h$
with $R$, let the roots in $R'$ be  
$\gamma_0,\gamma_1,\dots,\gamma_{2m-1}$ enumerated as in
\ref{dihedral}. The relations 
\begin{gather*}
     \sum_{i=0}^{m-1} [\gamma_i][\gamma_{i+k}] = 0
\quad\text{for all }k;
\\
  [\gamma_l]\cdot [\gamma_0][\gamma_1]\dots [\gamma_{2l}]
 + [\gamma_0][\gamma_1]\dots [\gamma_{2l}] \cdot [\gamma_l] 
+ [\gamma_l]\cdot [\gamma_{2l}][\gamma_{2l-1}]\dots [\gamma_0]
\\
+[\gamma_{2l}][\gamma_{2l-1}]\dots [\gamma_0] \cdot [\gamma_l] = 0
\quad\text{for }l=[m/2]-1,
\end{gather*}
are imposed in $\mathit{BE}(W,S)$. The second, $4$\dash term relation
is meaningful only when $m\ge 4$.

The bracket algebras generalise the quadratic algebras $\mathcal{E}_n$
to the case of arbitrary Coxeter group. 
If $W$ is a Weyl group of simply\dash laced type, one has
$\mathit{BE}(W,S) = \Nic_\qd(\VV{W})$ because there are no $4$\dash
term relations in $\mathit{BE}(W,S)$. 
The relations in the bracket algebra are just sufficient to prove that
$\mathit{BE}(W,S)$ 
contains a commutative subalgebra isomorphic to the coinvariant
algebra $S_W$, which was done in \cite{KM} for crystallographic
Coxeter groups of classical type and of type $G_2$. However, 
in some cases the bracket algebra  has `less' relations than 
$\Nic_W$ has. 

For example, when $W$ is the Weyl group of type $B_2$, one has $\dim
\mathit{BE}(W, S)=\infty$ according to \cite{KM}; $\dim
\Nic_{W} = 64$ which may be verified by a computer calculation.
In this case, $\Nic_W$ is the quotient of $\mathit{BE}(W,S)$
by the Coxeter relation $[\alpha_1][\alpha_2][\alpha_1][\alpha_2] = 
[\alpha_2][\alpha_1][\alpha_2][\alpha_1]$ between the simple root
generators. 

\begin{xproposition}
If $W$ is a Weyl group of type other than $G_2$,
the Nichols algebra $\Nic_W$ is a quotient of the bracket algebra
$\mathit{BE}(W,S)$.
\end{xproposition}
\begin{proof}
The relation (1) of the bracket algebra holds in $\Nic_W$.
Now let the root subsystem $R'=\{\gamma_0,\dots,\gamma_{2m-1}\}$ be as in (2).
One checks (see \ref{dihedral}) that $\Psi([\gamma_i]\otimes [\gamma_j])
=-[\gamma_{2i-j}]\otimes [\gamma_i]$ where $\Psi$ is the braiding
\ref{V_W} on $V_W$. 
Note
that $\gamma_{m+i}=-\gamma_i$ as the indices are taken modulo $2m$.
Applying  $\id+\Psi$ to the left hand side of the quadratic relation
in (2), one gets $\sum_{i=0}^m [\gamma_{i+k}]\otimes
[\gamma_i]-[\gamma_{i+2k}]\otimes [\gamma_{i+k}]$ which is zero. 
Thus, the quadratic relation in (2) holds in $\Nic_W$ 
because its left hand side
lies in $\ker(\id+\Psi)$. 

It remains to show that the 
$4$\dash term relation in (2) holds in $\Nic_W$. 
If $W$ is of type $A$, $B=C$, $D$, $E$ or $F$, a root
subsystem $R'\subset R$ of rank $2$ consists of at most $8$ roots,
i.e.\ $m\le 4$. For $m=4$ and $l=1$, the braided symmetriser
$[4]!_\Psi$, applied to the 
left hand side of the $4$\dash term relation, gives zero --- this is 
verified by easy computation using factorisation
\ref{braided_symmetriser} of $[4]!_\Psi$.
\end{proof}

It has been observed by T.~Maeno that in type $G_2$, the $4$\dash term
relation in the bracket algebra is not compatible with the braided
Hopf algebra structure and therefore cannot hold in the Nichols\dash
Woronowicz algebra. Thus, the statement of the Proposition is not true when $W$ is of type $G_2$.

\subsection{} 

The intriguing question remains, whether the Nichols algebra $\Nic_{\Symm_n}$
coincides with the quadratic algebra $\mathcal{E}_n$ or is a proper
quotient of it. 

The graded components of degrees $1$, $2$, $3$ in $\Nic_{\Symm_n}$ 
and $\mathcal{E}_n$ may be shown to coincide. 
Furthermore, $\Nic_{\Symm_n}=\mathcal{E}_n$ for $n\le 5$ (see
\cite[Example 6.4]{MS} for $n\le 4$, \cite{G} for $n=
5$). Incidentally, $\Symm_n$ for $n\le 5$ and $W_{B_2}$ are the only
examples of Coxeter groups where we know the Nichols algebra $\Nic_W$
to be finite\dash dimensional.

We finish with the following conjecture, which already appeared in a
number of sources including \cite{MS} and \cite{M}. If true, this
conjecture would mean that our construction of $\Nic_W$ as a model for
the Schubert calculus generalises, in proper sense, the Fomin\dash Kirillov
construction.

\begin{xconjecture} 
The algebras $\Nic_{\Symm_n}$ are quadratic and coincide with
the Fomin\dash Kirillov algebras $\mathcal{E}_n$ for all $n$.
\end{xconjecture}


\end{document}